\documentclass[11pt]{amsart}
\usepackage{amsmath,amscd}
\usepackage{amsfonts,amssymb}
\usepackage{latexsym}
\usepackage{amsbsy}
\usepackage[all]{xy}
\usepackage{amscd}
\usepackage[dvips]{graphicx}

\newtheorem{Prop}{Proposition}[section]
\newtheorem{Thm}[Prop]{Theorem}
\newtheorem{Lemma}[Prop]{Lemma}
\newtheorem{Cor}[Prop]{Corollary}

\newtheorem{Remark}[Prop]{Remark}

\newtheorem{definition}[Prop]{Definition}

\def\az{\alpha}      \def\ud{{\underline{d}}}

\def\lz{\lambda}

\def\bbn{{\mathbb N}}  \def\bbz{{\mathbb Z}}  \def\bbq{{\mathbb Q}} \def\bb1{{\mathbb 1}}
   \def\bbe{{\mathbb E}} 

  \def\bbc{{\mathbb C}}

\def\ra{\rightarrow}

\def\im{\mbox{Im}}

\def\hom{\mbox{Hom}}

\def\ext{\mbox{Ext}\,} 
\def\dim{\mbox{dim}\,}

\def\udim{\mbox{\underline {dim}}}
\def\ed{\mbox{End}}

\def\mod{\mbox{mod}\,}  
    
\def\Ker{\mbox{Ker}\,}

\def\Coker{\mbox{Coker}\,}

\def\uq2{U_q(\hat{sl}_2)}

\def\bb{{\bf b}}

\def\nd{{\noindent}}

\def\mf{{\mathcal{F}}}
\def\mt{{\mathcal{T}}}
\def\ms{{\mathcal{S}}}

\def\mp{{\mathcal{P}}}

\def\mo{{\mathcal{O}}}

\def\mt{{\mathcal{T}}}

\def\ue{{\underline{e}}}

\begin{document}

\title{Realizing Enveloping Algebras via Varieties of Modules}

\author{Ming Ding,  Jie Xiao and Fan Xu }
\address{Department of Mathematics, Tsinghua University,
 Beijing 100084, People's Repubic of China.}
\email{m-ding04@mails.tsinghua.edu.cn (M.Ding)\\ jxiao@math.tsinghua.edu.cn (J.Xiao)\\
f-xu04@mails.tsinghua.edu.cn (F.Xu)}

\dedicatory{Dedicated to Professor J.A.Green on the occasion of his
80th birthday }
\thanks{The research was supported in part by NSF of China and by
the 973 Project of the Ministry of Science and Technology of China.}
\date{}
\maketitle

\begin{abstract}
By using the Ringel-Hall algebra approach, we investigate the
structure of the Lie algebra $L(\Lambda)$ generated by
indecomposable constructible sets in the varieties of modules for
any finite dimensional $\mathbb{C}$-algebra $\Lambda.$ We obtain a
geometric realization of the universal enveloping algebra
$R(\Lambda)$ of $L(\Lambda).$ This generalizes the main result of
Riedtmann in \cite{R}. We also obtain Green's theorem in \cite{G}
in a geometric form for any finite dimensional
$\mathbb{C}$-algebra $\Lambda$ and use it to give the
comultiplication formula in $R(\Lambda).$
\end{abstract}

\section{Introduction}
Let us recall a  result of Riedtmann in \cite{R}. Let $\Lambda$ be a
finite dimensional associative $\mathbb{C}$-algebra with unit
element. If $\Lambda$  is of finite representation type, let
$\mathcal{J}$ be a set of representatives for the isomorphism
classes of indecomposable $\Lambda$-modules. Then
$\mathcal{P}=\{\bigoplus_{A\in\mathcal{J}}{A^{\mu(A)}}:\mu(A)\in
\mathbb{N}\}$ is a set of representatives for the isomorphism
classes of all $\Lambda$-modules. Let $R(\Lambda)$ be the free
$\mathbb{Z}$-module on the basis $\{u_{A},A\in \mathcal{P}\}$.
Define the multiplication  of two basis vectors by (see Section 2 in
\cite{R})
$$u_{A}\bullet u_{B}=\sum_{X\in\mathcal{P}
}{\chi(V(A,B;X))u_X},$$ where $\chi(V(A,B;X))$ is the
Euler-Poincar\'e characteristic of the variety:
$$V(A,B;X)=\{0\subseteq Z\subseteq X: Z\in \mathrm{mod} \Lambda,\ Z\cong A,
X/Z\cong B\}$$ This is a geometric version of Ringel-Hall algebras
(see \cite{Ri1}) by considering the Euler-Poincar\'e characteristic
of the filtration varieties to replace the filtration numbers over a
finite field (see also \cite{L} and \cite{S}). Then $R(\Lambda)$ is
an associative $\mathbb{Z}$-algebra with unit element. The
$\mathbb{Z}$-submodule $L(\Lambda)=\bigoplus_{A\in
\mathcal{J}}{\mathbb{Z}u_A}$ of $R(\Lambda)$ is a Lie subalgebra
with bracket $[x,y]=x\bullet y-y\bullet x$ (we will simply write
$xy$ for $x\bullet y$ in the below). She proved the following
result by a gradation method.

{\em For any numbering $\mathcal{J}=\{A_1,\cdots,A_n\},$  the
universal enveloping algebra $U(\Lambda)$ of $L(\Lambda)$, is freely
generated as a $\mathbb{Z}$-module by the classes of the words
$A_{1}^{\lambda_1}\cdots A_{n}^{\lambda_n},\ (\lambda_1,\cdots,
\lambda_n\in \mathbb{N})$. It is isomorphic to the subalgebra
$R(\Lambda)^\prime$ of $R(\Lambda)$ generated by $\{u_{A}:A\in
\mathcal{J}\}$, which is the free $\mathbb{Z}$-module generated by
the elements
$(\prod_{j=1}^{n}(\lambda_j{!}))u_{A_{1}^{\lambda_1}\oplus
\cdots\oplus A_{n}^{\lambda_n}}$}.

We note that a similar property for Ringel-Hall algebras
$H(\Lambda)$ has been obtained by Guo and Peng in \cite{GP} for any
finite dimensional algebra $\Lambda$ over a finite field.

In this paper we extend Riedtmann's result.  Let $R(\Lambda)$ be
the $\mathbb{Z}$-module generated by the characteristic functions
of constructible sets of stratified Krull-Schmidt (see Section 3.2
for the definition) and $L(\Lambda)$ be the $\mathbb{Z}$-submodule
of $R(\Lambda)$ generated by the characteristic functions of
indecomposable constructible sets (see Section 3.2 for the
definition). Our main result (Theorem \ref{Rie}) provides,
$R(\Lambda)$ can be realized as the enveloping algebra of
$L(\Lambda),$ from which Riedtmann's result can be deduced  by
replacing  orbits of modules by constructible sets.

Moreover, we propose a `degenerate form' (Theorem
\ref{greentheorem}) and an extended version (Remark
\ref{extended}) of Green's theorem on Hall algebras (\cite{G}, see
also \cite{Ri3}). The extended version holds for any finite
dimensional associative algebra (not only hereditary algebra). We
note that there is an analog of Green's theorem which was given by
Lusztig (\cite{L2}).

The paper is organized as follows. In Section 2 we recall the
basic knowledge about varieties of modules and constructible
functions. In Section 3 we give the definitions of $R(\Lambda)$
and $L(\Lambda)$ for any finite dimension algebra $\Lambda$ over
$\mathbb{C}$ which can be viewed as an extension of Section 2.3 in
\cite{R}. Following this, we study the relation between
$R(\Lambda)$ and $L(\Lambda)$ in Section 4. In particular, Theorem
4.2 is Riedtmann's result when considering representation-finite
case.  In Section 5, we give a formula in Theorem 5.7 which can be
viewed as a variant of Green's formula. In Section 6, we prove
that this formula induces the  comultiplication of $R(\Lambda)$ as
a Hopf algebra.

Finally, we note that there is a similar result in \cite{Joy} in
the context of constructible functions on stacks. The author use
quite different methods. It is interesting to consider analogues
of constructible subsets of stratified Krull-Schmidt for stacks.

\section{Basic Concepts}
In this section, we will fix notations and recall some elementary
concepts in algebraic geometry.

The set of all n-tuples with each coordinate element in the field
$\bbc$ is called an n-dimensional affine space over $\bbc$.  Each
closed set (with induced topology) of an affine space is called an
affine algebraic variety. A point $x$ of an affine variety $X$ is
smooth if the dimension of the tangent space of $X$ at $x$ is
equal to the dimension of $X$ at $x$, denoted by dim$T_{X,x}=$
dim$_{x}{X}$. Let $\varphi:X\longrightarrow Y$ be a dominating
(i.e $\varphi(X)$ is Zariski-dense in $Y$ ) regular map, if $x\in
X$ and $y=\varphi(x)\in Y$ is a smooth point, then we say that
$\varphi$ is smooth at $x$ if $x$ is a smooth point of $X$ and
d$\varphi$ maps $T_{x,X}$ onto $T_{y,Y}$. Let $X$ be a Zariski
topolgical space, the intersection of an open subset and a close
subset is called a locally closed subset. A subset in $X$ is
called constructible if it is a disjoint union of finitely many
locally closed subsets. Obviously, open sets and close sets are
both constructible sets.

\begin{definition}
A function $f$ on $X$ is called constructible if $X$ can be divided
into finitely many constructible sets satisfying that $f$ is
constant on each such constructible set.
\end{definition}

Let $\mathcal{O}$ be an above constructible set, $1_{\mathcal{O}}$
is called a characteristic function if $1_{\mathcal{O}}(x)=1$ for
any $x\in \mathcal{O}$; $1_{\mathcal{O}}(x)=0$ for any $x\notin
\mathcal{O}$. It is clear that $1_{\mathcal{O}}$ is the simplest
constructible function and any constructible function is a linear
combination of characteristic functions.

Let $\Lambda$ be a finite dimensional $\mathbb{C}$-algebra. By a
result of P.Gabriel in \cite{Ga}, the algebra $\Lambda$ is given by
a
 quiver $Q$ with relations $R$ (up to Morita equivalence).
Let $Q=(Q_0,Q_1,s,t)$ be a quiver, where $Q_0$ and $Q_1$ are the
sets of vertices and arrows respectively, and $s,t: Q_1\rightarrow
Q_0$ are
 maps such that any arrow $\alpha$ starts at $s(\alpha)$ and terminates at $t(\alpha).$
There are only finitely many simple $\Lambda$-modules
$S_1,\cdots,S_n$, up to isomorphism. The index set for simple
modules is denoted by $I$ (we may identify $I=Q_0$). For any
$\Lambda$-module $M$, we denote by $\udim M$ the vector in
$\mathbb{N}I$ whose $i$-th component is the multiplicity of $S_i$
in any composition series of $M$ for $i=1,\cdots,n$. It is called
the dimension vector of $M$. For any dimension vector
 $\underline{d}=\sum_i a_i i\in \mathbb{N}I,$ we consider the affine space over $\mathbb{C}$
$$\mathbb{E}_{\underline{d}}(Q)=\bigoplus_{\alpha\in Q_1}\hom_{\mathbb{C}}(\bbc^{a_{s(\az)}},\bbc^{a_{t(\az)}})$$
Any element $x=(x_{\az})_{\az\in Q_1}$ in $\bbe_{\ud}(Q)$ defines
a representation $M(x)=(\bbc^{\ud}, x)$ of $Q$ such that
$\bbc^{\ud}=\bigoplus_{i\in Q_0}\bbc^{a_i}$. A relation in $Q$ is
a linear combination $\sum_{i=1}^{r}\lz_{i}p_i,$ where
$\lz_i\in\bbc$ and $p_i$ are paths of length at least two with
$s(p_i)=s(p_j)$ and $t(p_i)=t(p_j)$ for all $1\leq i,j\leq r.$ For
any $x=(x_{\az})_{\az\in Q_1}\in\bbe_{\ud}(Q)$ and any path
$p=\az_1\az_2\cdots\az_m$ in $Q$, we set
$x_{p}=x_{\az_1}x_{\az_2}\cdots x_{\az_m}.$ Then $x$ satisfies a
relation
 $\sum_{i=1}^{r}\lz_{i}p_i$ if $\sum_{i=1}^{r}\lz_i x_{p_i}=0.$ If $R$ is a set of relations in $Q,$ then
let $\bbe_{\ud}(Q,R)$ be the closed subvariety of $\bbe_{\ud}(Q)$
which consists of elements satisfying all relations in $R.$ Any
element $x=(x_{\az})_{\az\in Q_1}$ in $\bbe_{\ud}(Q,R)$ defines in
a natural way a representation $M(x)$ of $\Lambda=\bbc Q/J$ with
$\udim M(x)=\ud,$ where $J$ is the admissible ideal generated by
$R.$
 We consider the algebraic group:
$$G_{\ud}:=G_{\ud}(Q)=\prod_{i\in I}GL(a_i,\bbc),$$ which acts on
 $\bbe_{\ud}(Q)$ by $(x_{\az})^{g}_{\az\in Q_1}=(g_{t(\az)}x_{\az}g_{s(\az)}^{-1})_{\az\in Q_1}$ for $g\in G_{\ud}$ and
  $(x_{\az})_{\az\in Q_1}\in\bbe_{\ud}(Q).$
It naturally induces the action of $G_{\ud}$ on $\bbe_{\ud}(Q,R).$
The orbit space is $\bbe_{\ud}(Q,R)/G_{\ud}.$ There is a natural
bijection between the set ${\mathcal M}(\Lambda, \bbc^{\ud})$ of
isomorphism classes of $\bbc$-representations of $\Lambda$ whose
underlying space is $\bbc^{\ud}$ and the set of orbits of
$G_{\ud}$ in $\bbe_{\ud}(Q,R).$ So we may identify  ${\mathcal
M}(\Lambda,\bbc^{\ud})$ with $\bbe_{\ud}(Q,R)/G_{\ud}.$ For any
constructible subset $\mathcal{O}$ in $\bbe_{\ud}(Q,R)$, we say
$\mathcal{O}$ is $G_{\ud}$-invariant if
$G_{\ud}\cdot\mathcal{O}=\mathcal{O}.$

In the following, the constructible sets and functions will always
be assumed $G_{\ud}$-invariant. For fixed $Q$, $R$ and $\ud$, we
simply write $\bbe_{\ud}$ instead of $\bbe_{\ud}(Q, R).$

Let $\mo$ and $\mo'$ be constructible subsets in $\bbe_{\ud}(Q,R)$
and $\bbe_{\ud'}(Q,R)$ respectively, we define
$$\mo\oplus\mo'=G_{\ud+\ud'}\cdot\{z\in\bbe_{\ud+\ud'}(Q,R)\mid
M(z)=M(x)\oplus M(x')\ \ \mbox{with}\ \ x\in\mo\ \mbox{and}\ x'$$
$$\in\mo'\}=\{z\in\bbe_{\ud+\ud'}(Q,R)\mid M(z)\cong M(x)\oplus M(x')\ \
\mbox{with}\ \ x\in\mo\ \mbox{and}\ x'\in\mo'\}.$$ We denote by
$n\mo$ the $n$ copies of direct sum of $\mo.$ Since
$$\psi:G_{\ud+\ud'}(Q)\times \bbe_{\ud}(Q,R)\times
\bbe_{\ud'}(Q,R)$$ $$\longrightarrow \bbe_{\ud}(Q,R)\oplus
\bbe_{\ud'}(Q,R)
$$
$$(g,V_1,V_2)\mapsto g\cdot (V_1\oplus V_2)$$ is an algebraic
map and also surjective, and the image of a constructible set is
constructible, we have that for any constructible subsets
$\mathcal{O},\mathcal{O'}$ in $\bbe_{\ud}(Q,R),\bbe_{\ud'}(Q,R)$
respectively, $\mo\oplus\mo'$ is a constructible subset in
$\bbe_{\ud+\ud'}(Q,R)$, that is to say, the direct sum of
constructible sets is a constructible set.

Let $\mathcal{O}_1,\cdots,\mathcal{O}_r$ be constructible sets and
$M$ be a $\Lambda$-module, we define the set:
$$V(\mathcal{O}_1,\cdots,\mathcal{O}_r;M)=\{0=M_0\subseteq M_1\subseteq\cdots \subseteq M_r=M\mid$$
$$M_i\in \mathrm{mod}\ \Lambda,M_i/M_{i-1}\in \mathcal{O}_i,1\leq i\leq r,i\in\mathbb{N}\}$$
In particular, $V(\bbe_{\ud_1},\bbe_{\ud_2};M)$ is the
Grassmannian of $\Lambda$-submodules of $M$, denoted by
$Gr_{\ud_1}(M).$

Let $\chi$ denote Euler characteristic in compactly-supported
cohomology. Let $X$ be an algebraic variety and $\mo$ a
constructible subset of $X$ as the disjoint union of finitely many
locally closed subsets $X_i$ for $i=1,\cdots,m.$ We have
$\chi(\mo)=\sum_{i=1}^m\chi(X_i).$  We will use the following
properties:
\begin{Prop}[\cite{R},\cite{Joyce}]\label{Euler} Let $X,Y$ be algebraic varieties over $\mathbb{C}.$
Then
\begin{enumerate}
    \item  If an algebraic variety $X$ is the disjoint union of
finitely many constructible sets $X_1,\cdots,X_r$, then
$$\chi(X)=\sum_{i=1}^{r}{\chi(X_i)}$$
    \item  If $\varphi:X\longrightarrow Y$ is a morphism
with the property that all fibers have the same Euler
characteristic $\chi$, then $\chi(X)=\chi\cdot \chi(Y).$ In
particular, if $\varphi$ is a locally trivial fibration in the
analytic topology with fibre $F,$ then $\chi(Z)=\chi(F)\cdot
\chi(Y).$
    \item $\chi(\bbc^n)=1$ and $\chi(\mathbb{P}^n)=n+1$ for all $n\geq
    0.$
\end{enumerate}
\end{Prop}
We recall {\it pushforward} functor from the category of algebraic
varieties over $\mathbb{C}$ and the category of
$\mathbb{Q}$-vector spaces ( see \cite{Mac} and \cite{Joyce}).

Let $\phi: X\rightarrow Y$ be a morphism of varieties. For $f\in
M(X)$ and $y\in Y,$ define
$$
\phi_{*}(f)(y)=\sum_{c\in\mathbb{Q}}c\chi(f^{-1}(c)\cap
\phi^{-1}(y))
$$
\begin{Thm}[\cite{Dimca},\cite{Joyce}]\label{Joyce}
Let $X,Y$ and $Z$ be algebraic varieties over $\mathbb{C},$ $\phi:
X\rightarrow Y$ and $\psi: Y\rightarrow Z$ be morphisms of
varieties, and $f\in M(X).$ Then $\phi_{*}(f)$ is constructible,
$\phi_{*}: M(X)\rightarrow M(Y)$ is a $\mathbb{Q}$-linear map and
$(\psi\circ \phi)_{*}=(\psi)_{*}\circ (\phi)_{*}$ as
$\mathbb{Q}$-linear maps from $M(X)$ to $M(Z).$
\end{Thm}

 In order to deal with orbit spaces, we also need to consider
the geometric quotients.
\begin{definition}\label{quotient}
Let $G$ be an algebraic group acting on a variety $X$ and
$\phi:X\rightarrow Y$ be a $G$-invariant morphism, i.e. a morphism
constant on orbits. The pair $(Y,\phi)$ is called a geometric
quotient if $\phi$ is open and for any open subset $U$ of $Y$, the
associated comorphism identifies the ring $\mo_{Y}(U)$ of regular
functions on $U$ with the ring $\mo_{X}(\phi^{-1}(U))^{G}$ of
$G$-invariant regular functions
 on $\phi^{-1}(U)$.
\end{definition}

The following result due to Rosenlicht \cite{Ro} is essential to
us.

\begin{Lemma}\label{Rosenlicht}
Let $X$ be a $G$-variety, then there exists an open and dense
$G$-stable subset which has a geometric $G$-quotient.
\end{Lemma}
By this Lemma, we can construct a finite stratification over $X.$
Let $U_1$ be an open and dense $G$-stable subset of $X$ as in
Lemma \ref{Rosenlicht}. Then
$dim_{\mathbb{C}}(X-U_1)<dim_{\mathbb{C}}X.$ We can use the above
lemma again, there exists a dense open $G$-stable subset $U_2$ of
$X-U_1$ which has a geometric $G$-quotient. Inductively, we get
the finite stratification $X=\cup_{i=1}^{l}U_{i}$ which $U_{i}$ is
$G$-invariant locally closed subset and has a geometric quotient,
$l\leq dim_{\mathbb{C}}X.$ We denote by $\phi_{U_i}$ the geometric
quotient map on $U_i.$ Define the \emph{quasi Euler-Poincar\'e
characteristic} of $X/G$ by
$$\chi(X/G):=\sum_{i}\chi(\phi_{U_i}(U_i)).$$ If $\{U'_j\}$ is
another choice for the definition of $\chi(X/G)$, then
$\chi(\phi_{U_i}(U_i))=\sum_{j}\chi(\phi_{U_i\cap U'_j}(U_i\cap
U'_j))$ and $\chi(\phi_{U'_j}(U'_j))=\sum_{i}\chi(\phi_{U_i\cap
U'_j}(U_i\cap U'_j)).$ Thus
$$\sum_{i}\chi(\phi_{U_i}(U_i))=\sum_{j}\chi(\phi_{U'_j}(U'_j)).$$
Hence, $\chi(X/G)$ is well-defined (see \cite{XXZ}). Similarly, we
obtain $\chi(\mo/G):=\sum_i\chi(\phi_{U_i}(\mo\bigcap U_i))$ is
well-defined for any $G$-invariant constructible subset $\mo$ of
$X.$

We also introduce the following notation. Let $f$ be a
constructible function over a variety $X,$ it is natural to define
\begin{equation}\label{integral}
\int_{x\in X}f(x):=\sum_{m\in \bbc}m\chi(f^{-1}(m))
\end{equation}
Comparing with Proposition \ref{Euler}, we also have
\begin{Prop}[\cite{XXZ}]\label{Euler2} Let $X,Y$ be algebraic varieties over
$\mathbb{C}$ under the actions of the algebraic groups $G$ and $H$
respectively.  Then
\begin{enumerate}
    \item  If an algebraic variety $X$ is the disjoint union of
finitely many $G$-invariant constructible sets $X_1,\cdots,X_r$,
then
$$\chi(X/G)=\sum_{i=1}^{r}{\chi(X_i/G)}$$
    \item  If $\varphi:X\longrightarrow Y$ is a morphism induces
    the quotient map $\phi:X/G\rightarrow Y/H$
with the property that all fibers for $\phi$ have the same Euler
characteristic $\chi$, then $\chi(X/G)=\chi\cdot \chi(Y/H).$
   \end{enumerate}
\end{Prop}

Moreover, if there exists an action of algebraic group $G$ on $X$
as in Definition \ref{quotient}, and $f$ is $G$-invariant
constructible function over $X,$ We define
\begin{equation}\label{integral2}
\int_{x\in X/G}f(x):=\sum_{m\in \bbc}m\chi(f^{-1}(m)/G)
\end{equation}

\section{Definitions of $R(\Lambda)$ and $L(\Lambda)$}

\nd {\bf 3.1} \ \ Let $\mathcal{O}_i$ be a $G_{\ud_i}$-invariant
constructible subset of  $\bbe_{\underline{d}_i}$ for i=1,2 and
$\mathcal{L}$ be a constructible subset  of $\bbe_{\ud_1+\ud_2}$
(need not be invariant under the group action). Similar to
$V(\mathcal{O}_1,\mathcal{O}_2;M)$, we define:
$$V(\mathcal{O}_1,\mathcal{O}_2;\mathcal{L})=\{(L,L_1)\mid L\in \mathcal{L}, L_1\in V(\mo_1,\mo_2;L)\}. $$ There exists a natural morphism of
varieties:
$$
\widetilde{\pi}:
V(\mathcal{O}_1,\mathcal{O}_2;\mathcal{L})\rightarrow \mathcal{L}
$$
mapping $(L, L_1)$ to $L.$ It is clear to see
$\chi(\widetilde{\pi}^{-1}(L))=\chi(V(\mathcal{O}_1,\mathcal{O}_2;L))$.

\begin{definition} $1_{\mathcal{O}_1}\bullet
1_{\mathcal{O}_2}(L)=\chi(V(\mathcal{O}_1,\mathcal{O}_2;L))$
\end{definition}

Consider the $\bbc$-space
$M_G(\Lambda)=\bigoplus_{\ud\in\bbn^n}M_{G_{\ud}}(Q,R)$ where
$M_{G_{\ud}}(Q,R)$ is the $\bbc$-space of $G_{\ud}$-invariant
constructible function on $\bbe_{\ud}(Q,R).$ Actually we have
defined the convolution multiplication on $M_G(\Lambda)$ as follows.
For any $f\in M_{G_{\ud}}(Q,R)$ and $g\in M_{G_{\ud'}}(Q,R),$
$f\bullet g\in M_{G_{\ud+\ud'}}(Q,R)$ is given by the formula
$$f\bullet
g(L)=\sum_{c,d\in\bbc}\chi(V(f^{-1}(c),g^{-1}(d);L))cd.$$ As usual
for an algebraic variety $V$ and a constructible function $f$ on
$V,$ we define $\int_{V}f=\sum_{c\in\bbc}\chi(f^{-1}(c))c.$ Define
$1_{0}$ be the characteristic function satisfying $1_0(X)=1$ if
$X\cong 0$ and $1_0(X)=0$ otherwise. Similar as in \cite{R} (see
also Section 2 in \cite{L3}), we can prove that:

\begin{Prop}\label{associativity} The space $M_G(\Lambda)$ under the convolution
multiplication $\bullet$ is an associative algebra over $\bbc$
with identity $1_0$.
\end{Prop}
\begin{proof}
It is enough to prove that $1_{\mo_1}\bullet 1_{\mo_2}$ is
constructible for any constructible subsets $\mo_1$ and $\mo_2$ of
$\bbe_{\ud_1}$ and $\bbe_{\ud_2},$ respectively.
 We consider the canonical morphism of varieties
$$
\widetilde{\pi}: V(\mo_1, \mo_2 ;\bbe_{\ud_1+\ud_2})\rightarrow
\bbe_{\ud_1+\ud_2}.
$$
We have
$\widetilde{\pi}^{-1}(L)=V(\mathcal{O}_1,\mathcal{O}_2;L).$ By
Theorem \ref{Joyce}, we obtain a constructible function
$$\widetilde{\pi}_{*}(1_{V(\mo_1, \mo_2
;\bbe_{\ud_1+\ud_2})})=1_{\mo_1}\bullet 1_{\mo_2}.$$ For any $L\in
\bbe_{\ud_1+\ud_2}$ and $\mo_i\subset \bbe_{\ud_i}$ for $i=1,2,3,$
define
$$
\mathcal{A}_1=\{L_1\subseteq L_2\subseteq L\mid L_1\in \mo_2,
L_2/L_1\in \mo_1, L/L_2\in \mo_3\}
$$
and
$$
\mathcal{A}_2=\{( L_1\subseteq L,L_2/L_1\subseteq L/L_1)\mid
L_1\in \mo_2, L_2/L_1\in \mo_1, L/L_2\in \mo_3\}.
$$
Then
$$
\chi(\mathcal{A}_1)=\chi(\mathcal{A}_2).
$$
Furthermore, we have
$$(1_{\mo_2}\bullet 1_{\mo_1})\bullet 1_{\mo_3}(L)=\chi(\mathcal{A}_1)$$
and
$$
1_{\mo_2}\bullet (1_{\mo_1}\bullet
1_{\mo_3})(L)=\chi(\mathcal{A}_2).
$$
Hence,
$$
1_{\mo_2}\bullet(1_{\mo_1}\bullet 1_{\mo_3})=(1_{\mo_2}\bullet
1_{\mo_1})\bullet 1_{\mo_3}.
$$
\end{proof}

\bigskip

\nd {\bf 3.2} \ \  A constructible set is called indecomposable if
all points in it correspond to indecomposable
$\Lambda$-modules(\cite{XXZ}). Let $\mathcal{O}$ be a
constructible set. If it has  the form:
$\mo=n_{1}\mathcal{O}_{1}\oplus \cdots \oplus
n_{k}\mathcal{O}_{k}$, where $\mathcal{O}_{i},1\leq i\leq k$ are
indecomposable constructible sets, then $\mo$ is called  to be of
Krull-Schmidt. A constructible set $\mathcal{Q}$ is called to be
of stratified Krull-Schmidt if it has a finite stratification
$\mathcal{Q}=\dot{\bigcup}_{i}\mathcal{Q}_i$ where each
$\mathcal{Q}_i$ is locally closed in $\mathcal{Q}$ and is of
Krull-Schmidt. Define
$R(\Lambda)=\Sigma{\mathbb{Z}1_{\mathcal{Q}}}$, where
$\mathcal{Q}$ are  constructible sets of stratified Krull-Schmidt.
Obviously $R(\Lambda)=\Sigma{\mathbb{Z}1_{\mathcal{O}}}$, where
$\mathcal{O}$ are  constructible sets of  Krull-Schmidt. Note that
for a Krull-Schmidt constructible set
$\mo=n_{1}\mathcal{O}_{1}\oplus \cdots \oplus
n_{k}\mathcal{O}_{k}$, we may assume $\mathcal{O}_{i},1\leq i\leq
k$ are disjoint to each other. For if ${\mathcal{O}_i\cap
\mathcal{O}_j}\neq \emptyset$, then
$$\mathcal{O}_i\oplus \mathcal{O}_j=
2(\mathcal{O}_i\cap
\mathcal{O}_j)\dot{\bigcup}\big((\mathcal{O}_i\backslash(\mathcal{O}_i\cap
\mathcal{O}_j))\oplus (\mathcal{O}_j\backslash(\mathcal{O}_i\cap
\mathcal{O}_j)))$$
$$\dot{\bigcup}\big((\mathcal{O}_i\cap \mathcal{O}_j)
\oplus(\mathcal{O}_j\backslash(\mathcal{O}_i\cap
\mathcal{O}_j)))\dot{\bigcup}\big((\mathcal{O}_i\cap
\mathcal{O}_j)\oplus (\mathcal{O}_i\backslash(\mathcal{O}_i\cap
\mathcal{O}_j)))$$ where $\dot{\bigcup}$ are  disjoint unions.

\begin{Lemma}\label{nonzero}
Let $\mo_1,\mo_2$ be two indecomposable constructible subsets of
$\bbe_{\ud_1}$ and $\bbe_{\ud_2}$. Let $L$ be a $\Lambda$-module
with dimension vector $\ud_1+\ud_2$. If
$\chi(V(\mo_1,\mo_2;L))\neq 0$, then there exist $X\in \mo_1, Y\in
\mo_2$ satisfying $$\chi(V(X,Y;L)) \neq 0.$$
\end{Lemma}
\begin{proof} For any submodule $L_1$ with dimension vector $\ud_1$ of $L$, by the knowledge of linear algebra, there exist unique $(\bbc^{\ud_1},x)\in \bbe_{\ud_1}$  isomorphic to $L_1$
and $(\bbc^{\ud_2},x')\in \bbe_{\ud_2}$ isomorphic to $L/L_1.$ In
fact, it induces a morphism of varieties (see also Lemma
\ref{kernel})
$$ \xymatrix {Gr_{\ud_1}(L)\ar[r]^-{\pi_1} &\bbe_{\ud_1}\times
\bbe_{\ud_2}}.
$$

For any constructible subsets $\mo_1$ and $\mo_2$ of
$\bbe_{\ud_1}$ and $\bbe_{\ud_2},$ respectively, we consider the
morphism:
$$
\xymatrix {Gr_{\ud_1}(L)\ar[r]^-{\pi_1} &\bbe_{\ud_1}\times
\bbe_{\ud_2}\ar[r]^-{\pi_2}& \bigcup_{i}\phi_i(U_i)}
$$
where $\bbe_{\ud_1}\times \bbe_{\ud_2}=\bigcup_{i}U_i$ is a finite
stratification with respect to algebraic group $G_{\ud_1}\times
G_{\ud_2}$ actions and $\phi_i: U_i\rightarrow \phi_i(U_i)$ is the
geometric quotient for any $i,$ and  $\pi_2=\bigcup_i\phi_i.$ It
is clear that
$$
(\pi_2\pi_1)^{-1}(\pi_2(\mo_1\times\mo_2))=V(\mo_1,\mo_2;L)
$$
By the proof of Proposition \ref{associativity},
$\{\chi(V(X,Y;L))\mid X\in \mo_1, Y\in \mo_2\}$ is a finite set.
Using Proposition \ref{Euler2} for the composition of morphisms
$\pi_2\pi_1,$ there exists a partition
$\pi_2(\mo_1\times\mo_2)=\bigcup \mp_i$ satisfying that for fixed
$\mp_i,$ $\chi(V(X,Y;L))$ is constant for any $(X,Y)\in \mp_i$
(denoted by $\chi_i$) and
$$
\chi(V(\mo_1,\mo_2;L))=\sum_i\chi(\mp_i)\cdot \chi_i\neq 0
$$
Hence, there exists some $\chi_i\neq 0.$ Thus we completed our
proof.
\end{proof}

\begin{Lemma}\label{zero}
Let $\mathcal{O}_1,\mathcal{O}_2$ be indecomposable constructible
sets as above. Let $M$ and $N$ be two nonzero $\Lambda$-modules
such that one of them does not belong to $\mo_1\cup\mo_2$. Then
$$\chi(V(\mathcal{O}_1,\mathcal{O}_2;M\oplus N))=0.$$
\end{Lemma}
\begin{proof} Let $L=M\oplus N$. If $\chi(V(\mathcal{O}_1,\mathcal{O}_2;L))\neq 0$, then from
Lemma \ref{nonzero}, there exist $A\in \mathcal{O}_1$ and $B\in
\mathcal{O}_2 $ such that $\chi(V(A,B;L))\neq 0$. From Lemma 2.2
in \cite{R}, we know $L\cong A\oplus B$. However according to
Krull-Schmidt theorem: $A\cong M,B\cong N$; or $A\cong N,B\cong
M$. It is a contradiction.
\end{proof}

\begin{Remark}
If $\mathcal{O}_1\cap \mathcal{O}_2=\emptyset$, then according to
Lemma \ref{zero}, we have
$$1_{\mathcal{O}_1}\bullet
1_{\mathcal{O}_2}=1_{\mathcal{O}_{1}\oplus
\mathcal{O}_2}+\sum_{i=1}^m c_i 1_{\mathcal{P}_i}$$ where
$\mathcal{P}_i$ are indecomposable constructible sets and
$c_i=\chi(V(\mathcal{O}_1,\mathcal{O}_2;L_i))$ for any
$L_i\in\mp_i.$
\end{Remark}

For $M\in\mod\Lambda,$ let $\gamma(M)$ be the number of
indecomposable direct summands of the Krull-Schmidt decomposition
of $M.$ For a constructible set $\mo$ we define
$\gamma(\mo)=max\{\gamma(M)\mid M\in\mo\}.$

According to Lemma \ref{zero} and Lemma 2.2 in \cite{R}, for
indecomposable constructible sets $\mo_1,\cdots,\mo_k$ and a
$\Lambda$-module $L,$
$\chi(V(\mathcal{O}_1,\cdots,\mathcal{O}_k;L))\neq 0$ implies that
$\gamma(L)\leq k$ and if $\gamma(L)=k$ then
$L=M_1\oplus\cdots\oplus M_k$ with $M_i\in\mo_i,\ i=1,\cdots,k.$

\begin{Prop} Let $\mathcal{O}$ be an indecomposable constructible set. Then $1_{\mathcal{O}}^k=
k!1_{k\mathcal{O}}+\sum_i{m_{i}1_{\mathcal{P}_i}}$, where
$\gamma(\mathcal{P}_i)<k.$
\end{Prop}
\begin{proof}
It is obvious for $k=1.$ If $k=2$, then
$1_\mathcal{O}^2(A\oplus A)=\chi(V(\mo,\mo;A\oplus
A))=\chi(V(A,A;A\oplus A))=2$ and $1_{\mo}^2(A\oplus
B)=\chi(V(\mo,\mo;A\oplus B))=\chi(V(A,B;A\oplus
B))+\chi(V(B,A;A\oplus B))=2$ for $A,B\in\mo$ with  $A\not\cong
B.$ Then $1_{\mo}^2=2\cdot1_{\mo\oplus\mo}+\sum_i m_i 1_{\mp_i}$
where $\mp_i$ are indecomposable constructible sets.

Suppose that the proposition is true for $k=n$, then for $k=n+1$, we
have $1_{\mathcal{O}}^{n+1}=1_{\mathcal{O}}^{n}1_\mathcal{O}$. Our
aim is to show the initial monomial of $1_{\mathcal{O} }^{n+1}$ is
$(n+1)!1_{(n+1)\mathcal{O}}$, so it is sufficient for us to prove
the initial monomial of $n!1_{n\mathcal{O}}1_\mathcal{O}$ is
$(n+1)!1_{(n+1)\mathcal{O}}$. Namely the initial monomial of
$1_{n\mo}1_{\mo}$ is $(n+1)1_{(n+1)\mo}.$ Consider any numbering
$\{A,B,\cdots,L\}$ in $\mo$  which  are not isomorphic to each
other, and natural numbers $a,b,\cdots,l$ satisfying
$a+b+\cdots+l=n$. Then $aA\oplus bB\oplus\cdots\oplus lL\in n\mo.$
Since
$$1_{aA\oplus bB\oplus\cdots\oplus
lL}1_A=(a+1)1_{(a+1)A\oplus bB\oplus\cdots\oplus lL}+\cdots$$
$$1_{(a+1)A\oplus
(b-1)B\oplus\cdots\oplus lL}1_B=b1_{(a+1)A\oplus
bB\oplus\cdots\oplus lL}+\cdots$$
$$\vdots$$
$$1_{(a+1)A\oplus
bB\oplus\cdots\oplus (l-1)L}1_L=l1_{(a+1)A\oplus
bB\oplus\cdots\oplus lL}+\cdots.$$ By the sum of the above
equalities, we obtain that
$$1_{n\mo}1_{\mo}((a+1)A\oplus bB\oplus\cdots\oplus
lL)=(a+1)+b+\cdots+l=n+1.$$  Thus we finished the proof.
\end{proof}

\begin{Cor}\label{coro3.7} Let $\mathcal{O}_1,\cdots,\mathcal{O}_k$ be indecomposable
constructible sets  and they are disjoint to each other. Then the
initial term of $1_{\mathcal{O}_1}^{n_1}\cdots
1_{\mathcal{O}_k}^{n_k}$ is
$$n_1!\cdots n_k!1_{n_1\mathcal{O}_1\oplus\cdots\oplus
n_k\mathcal{O}_k},$$ and the initial term of
$1_{m_1\mo_1\oplus\cdots\oplus m_k\mo_k}
1_{n_1\mo_1\oplus\cdots\oplus n_k\mo_k}$ is
$$\prod_{i=1}^k\frac{(m_i+n_i)!}{m_i!n_i!} 1_{(m_1+n_1)\mo_1\oplus\cdots\oplus (m_k+n_k)\mo_k},$$

\end{Cor}

\begin{Prop}\label{action}
 Let
$\mathcal{O}_1,\mathcal{O}_2$ be two constructible sets  of
Krull-Schmidt. Then we have the finite sum

$$1_{\mathcal{O}_1}\bullet 1_{\mathcal{O}_2}
=\sum_{\gamma(\mathcal{Q}_i)\leq
\gamma(\mo_1)+\gamma(\mo_2)}{q_i1_{\mathcal{Q}_i}}
$$ where  $\mathcal{Q}_i$ are constructible sets of stratified Krull-Schmidt.
\end{Prop}
\begin{proof}
As we know, $1_{\mathcal{O}_1}\bullet 1_{\mathcal{O}_2}$ is a
constructible function, so
$$
1_{\mathcal{O}_1}\bullet 1_{\mathcal{O}_2}=\sum_i
q_i1_{\mathcal{Q}_i}
$$
for some constructible sets $\mathcal{Q}_i$ where
$q_i=\chi(V(\mo_1,\mo_2;L_i))\neq 0$ for any $L_i\in
\mathcal{Q}_i$. By Lemma \ref{nonzero}, there exist $ X_i\in
\mo_1,Y_i\in \mo_2$ such that $\chi(V(X_i,Y_i;L_i))\neq 0$, then
$\gamma(L_i)\leq \gamma(X_i)+\gamma(Y_i)$ and  $L_i\cong X_i\oplus
Y_i$ if $\gamma(L_i)=\gamma(X_i)+\gamma(Y_i)$ by Lemma 2.2 in
\cite{R}, which implies $\gamma(\mathcal{Q}_i)\leq
\gamma(\mo_1)+\gamma(\mo_2).$

Moreover, for fixed $L_i$, suppose $L_i=\oplus_{j=1}^{s}Z_{j}$ for
some indecomposables $Z_1,\cdots,Z_s$. Similar to Lemma 2.2 in
\cite{R}, we define an action of $\mathbb{C}^{*}$ on $L_i$ by
$$
t\cdot (z_1,\cdots,z_s)=(tz_1,\cdots,t^{s}z_s)
$$
for $t\in \mathbb{C}^{*}$ and $z_j\in Z_j$, $j=1\cdots s.$ A
subspace $L'_i\subseteq L_i$ is stable under this action if and only
if $L'_i=\oplus_{j=1}^{s}(L'_i\cap Z_j).$ This induces the action of
$\mathbb{C}^{*}$ on $V(X_i,Y_i;L_i)$ whose fixed point is filtration
as follows:
$$0=
L_{i0}\subseteq L_{i1}\subseteq L_{i2}=L_i
$$
where $L_{i1}=\oplus_{j} (L_{i1}\cap Z_j)\cong X_i$ and
$L_i/L_{i1}=\oplus_{j}Z_j/(L_{i1}\cap Z_j)\cong Y_i$. The fixed
point set is not empty just because $V(X_i,Y_i;L_i)$ has the same
Euler-Poincar\'e characteristic as its fixed point set under the
action of $\mathbb{C}^{*}$ and we have assumed that
$\chi(V(X_i,Y_i;L_i)))\neq 0.$ We have the following exact
sequence:
$$
0\longrightarrow L_{i1}\cap Z_j\longrightarrow Z_j\longrightarrow
Z_j/(L_{i1}\cap Z_j)\longrightarrow 0
$$
for $j=1,\cdots,s.$ This means all indecomposable direct summands of
$L_i$ are the extension of direct summands of $X_i$ and $Y_i.$

Without loss of generality, we may assume that
$$\mo_1=\bigoplus_{i=1}^t m_i\mp_i \ \  \mbox{and}\ \
\mo_2=\bigoplus_{i=1}^t n_i\mp_i$$ with the property that each
$\mp_i$ is indecomposable and $\mp_i\cap\mp_j=\emptyset$ or
$\mp_i=\mp_j$ for all $i\neq j,$ and $m_i, n_j$ equal $0$ or $1$ for
all $i$ and $j.$ We call $\{A_1,\cdots,A_s\}$ an $s$-partition of
$\{1,2,\cdots,t\}$ if $A_1\cup\cdots\cup A_s=\{1,2,\cdots,t\}$,
$A_i\neq\emptyset$ for all $i$ and $A_i\cap A_j=\emptyset$ for all
$i\neq j$. Obviously the number of all partitions of
$\{1,2,\cdots,t\}$ is finite. For any two $s$-partitions
$\{A_1,\cdots,A_s\}$, $\{B_1,\cdots,B_s\}$ of $\{1,2,\cdots,t\}$ and
$z_l\in\bbz$, $l=1,\cdots,s,$ set
$$\ms_{A_l,B_l,z_l}=\{M\ \ \mbox{indecomposable}\ \mid
\chi(V(\bigoplus_{i\in A_l}m_i\mp_i,\bigoplus_{j\in B_l}n_j\mp_j;
M))=z_l\}$$
$$\mt_{A_l,B_l,z_l}=\{M\mid M\in \bigoplus_{i\in A_l}m_i\mp_i\oplus\bigoplus_{j\in
B_l}n_j\mp_j,$$ $$ \chi(V(\bigoplus_{i\in
A_l}m_i\mp_i,\bigoplus_{j\in B_l}n_j\mp_j; M))=z_l\}$$ These are
constructible sets, and by Corollary \ref{coro3.7}, are of
Krull-Schmidt. According to our analysis of $\bbc^*$ action, we
see that $1_{\mo_1}\bullet 1_{\mo_2}$ is a finite $\bbz$-linear
sum of $1_{\oplus_{l=1}^s\mo_{A_l,B_l,z_l}},$ where
$\mo_{A_l,B_l,z_l}=\ms_{A_l,B_l,z_l}$ or $\mt_{A_l,B_l,z_l}$ for
all $s$-partitions $(s=1,\cdots,t).$ The proof is finished.
\end{proof}

\begin{Thm} The $\bbz$-module $R(\Lambda)$ under the convolution multiplication $\bullet$  is an associative
$\mathbb{Z}$-algebra with unit element.
\end{Thm}

\nd{\bf 3.3} \ \ Define the $\mathbb{Z}$-submodule
$L(\Lambda)=\sum_{i}{\mathbb{Z}1_{\mathcal{O}_i}}$ of
$R(\Lambda)$, where $\mathcal{O}_i$ are indecomposable
constructible sets. Then we have the following result.

\begin{Thm} The $\bbz$-submodule $L(\Lambda)$ is a Lie subalgebra of $R(\Lambda)$ with bracket
$[x,y]=xy-yx$.
\end{Thm}
\begin{proof} We need to verify that $L(\Lambda)$ is closed under the Lie
bracket. If $\mathcal{O}_1=\mathcal{O}_2$, then
$1_{\mathcal{O}_1}\bullet
1_{\mathcal{O}_2}-1_{\mathcal{O}_2}\bullet 1_{\mathcal{O}_1}=0\in
L(\Lambda)$

According to Subsection 3.2, we may assume that $\mathcal{O}_1\cap
\mathcal{O}_2=\emptyset.$  By Remark 3.5 we have
$$1_{\mathcal{O}_1}\bullet
1_{\mathcal{O}_2}=1_{\mathcal{O}_{1}\oplus
\mathcal{O}_2}+\sum_{i}{\chi(V(\mathcal{O}_1,\mathcal{O}_2;L_i))1_{\mathcal{P}_i}}$$
$$1_{\mathcal{O}_2}\bullet
1_{\mathcal{O}_1}=1_{\mathcal{O}_{1}\oplus
\mathcal{O}_2}+\sum_{j}{\chi(V(\mathcal{O}_2,\mathcal{O}_1;N_j))1_{\mathcal{P}_j^{\prime}}}$$
where $\mathcal{P}_i,\mathcal{P}_{j}^{\prime}$ are indecomposable
constructible sets. So $1_{\mathcal{O}_1}\bullet
1_{\mathcal{O}_2}-1_{\mathcal{O}_2}\bullet 1_{\mathcal{O}_1}\in
L(\Lambda).$
\end{proof}

Note that $1_{\bbe_{\ud}(Q,R)}$ always lies in $R(\Lambda),$ by
the following lemma.
\begin{Lemma}\label{stratified KS} The variety $\bbe_{\ud}(Q,R)$ is of stratified
Krull-Schmidt for any $\ud\in\bbn I.$
\end{Lemma}

For any $\Lambda$-module $M$ with dimension vector $\ud$, we have
the Krull-Schmidt decomposition $M\cong\bigoplus_{i=1}^s m_iM_i.$
Let $\ud_i=\udim M_i.$ Then we call the vectors
$\{m_i\ud_i\}_{i=1}^s$ the Krull-Schmidt vectors of $M,$ also of
$\ud.$ Note that, for a module $M$ its Krull-Schmidt vectors are
uniquely determined and the dimension vector $\ud$ only has finitely
many families of Krull-Schmidt vectors. This notion coincides with
the Konstant partition if we consider quivers without relations.

\nd{\it Proof of Lemma \ref{stratified KS}}: Define
$$\bbe^{ind}_{\ud_i}=\{M\in\bbe_{\ud_i}(Q,R)\mid M\ \ \mbox{is
indecomposable}\}.$$ It is a locally closed subset. Then it is easy
to see that, by Krull-Schmidt theorem,
$$\bbe_{\ud}(Q,R)=\bigcup_{\{m_i\ud_i\}_{i=1}^s}\bigoplus_{i=1}^s
m_i\bbe^{ind}_{\ud_i}.$$ The union is a disjoint union and runs over
all Krull-Schmidt vectors of $\ud.$ The proof is finished. $\Box$

\section{The Universal Enveloping Algebra of $L(\Lambda)$}

 We consider the following tensor algebra over
$\bbz:$
$$T(L(\Lambda))=\bigoplus_{i=0}^{\infty}L(\Lambda)^{\otimes i}$$
where $L(\Lambda)^{\otimes 0}=\bbz,$ $L(\Lambda)^{\otimes
i}=L(\Lambda)\otimes\cdots\otimes L(\Lambda)$ ( $i$ times ). Then
$T(L(\Lambda))$ is an associative $\bbz$-algebra using the tensor as
multiplication. Let $J$ be the two-sided ideal of $T(L(\Lambda))$
generated by
$$1_{\mo_1}\otimes 1_{\mo_2}-1_{\mo_2}\otimes 1_{\mo_1}-[1_{\mo_1},1_{\mo_2}]$$
where $\mo_1,\ \mo_2$ are any indecomposable constructible sets
and $[-,-]$ is Lie bracket in $L(\Lambda).$ Then
$U(\Lambda)=T(L(\Lambda))/J$ is the universal enveloping algebra
of $L(\Lambda)$ over $\bbz.$ To avoid confusion, we write the
multiplication in $U(\Lambda)$ as $\ast$ and the multiplication in
$R(\Lambda)$ as $\bullet.$ We have the canonical homomorphism
$\varphi:U(\Lambda)\rightarrow R(\Lambda)$ satisfying
$\varphi(1_{\mo_1}\ast 1_{\mo_2})=1_{\mo_1}\bullet 1_{\mo_2}.$

\begin{Lemma} The canonical homomorphism
$\varphi:U(\Lambda)\rightarrow R(\Lambda)$ is an embedding.
\end{Lemma}
\begin{proof} We arbitrarily choose indecomposable constructible sets
$\mo_1,\cdots,\mo_n$ such that they are disjoint to each other. Then
$1_{\mo_1}, \cdots, 1_{\mo_n}$ are $\bbz$-linear independent in
$L(\Lambda).$ Set
$$U(\Lambda)_{\mo_1\cdots \mo_n}=\{\sum\lambda_w w\in
U(\Lambda)\mid w=1_{\mo_1}^{\ast\lambda_1}\ast\cdots\ast
1_{\mo_n}^{\ast\lambda_n}, $$$$\lambda_i\geq 0,\ \mbox{for}\
i=1,\cdots,n\ \mbox{and}\ \lambda_w\in\bbz\}$$ and let
$R(\Lambda)_{\mo_1\cdots\mo_n}$ be the subalgebra of $R(\Lambda)$
generated by $1_{\lambda_1\mo_1\oplus\cdots\oplus \lambda_n\mo_n}$
for all $\lambda_i\geq 0, i=1,\cdots,n.$ Then we have the
restriction map
$$\varphi\mid_{U(\Lambda)_{\mo_1\cdots\mo_n}}:
U(\Lambda)_{\mo_1\cdots\mo_n}\rightarrow
R(\Lambda)_{\mo_1\cdots\mo_n}$$ We will prove that
$\varphi\mid_{U(\Lambda)_{\mo_1\cdots\mo_n}}$ is an embedding. Set
$$U(\Lambda)_{\mo_1\cdots \mo_n}^{(m)}=\{\sum\lambda_w w\in
U(\Lambda)\mid w=1_{\mo_1}^{\ast\lambda_1}\ast\cdots\ast
1_{\mo_n}^{\ast\lambda_n} \ \mbox{with}\ \sum_{i=1}^n\lambda_i\leq
m\}$$ for any $m\in\bbn.$ By the PBW theorem, we know that
$$U(\Lambda)_{\mo_1\cdots \mo_n}^{(m)}/U(\Lambda)_{\mo_1\cdots
\mo_n}^{(m-1)}$$ is a free $\bbz$-module with the basis
$$
\{1_{\mo_1}^{\ast\lambda_1}\ast\cdots\ast
1_{\mo_n}^{\ast\lambda_n}\mid\sum_{i=1}^n\lambda_i=m\}.$$ Set
$$R(\Lambda)_{\mo_1\cdots\mo_n}^{(m)}=\{\sum \lambda_{\mo}1_{\mo}\mid
1_{\mo}\in R(\Lambda)_{\mo_1\cdots\mo_n}, \lambda_{\mo}\in\bbz \
\mbox{and}\ \gamma(\mo)\leq m\}$$ for any $m\in\bbn.$ It is known
that, by Krull-Schmidt theorem,
$$\{1_{\lambda_1\mo_1\oplus\cdots\oplus\lambda_n\mo_n}\mid
\lambda_i\geq 0 \ \mbox{for}\ i=1,\cdots,n\ \mbox{and}\
\sum_{i=1}^n\lambda_i=m\}$$ are $\bbz$-linear independent in
$$R(\Lambda)_{\mo_1\cdots \mo_n}^{(m)}/R(\Lambda)_{\mo_1\cdots
\mo_n}^{(m-1)}.$$ By Corollary \ref{coro3.7} the induced map of
$\varphi$
$$\bar{\varphi}:U(\Lambda)_{\mo_1\cdots \mo_n}^{(m)}/U(\Lambda)_{\mo_1\cdots
\mo_n}^{(m-1)}\rightarrow R(\Lambda)_{\mo_1\cdots
\mo_n}^{(m)}/R(\Lambda)_{\mo_1\cdots \mo_n}^{(m-1)}$$ sends
$1_{\mo_1}^{\ast\lambda_1}\ast\cdots\ast 1_{\mo_n}^{\ast\lambda_n}$
to
$\lambda_1!\cdots\lambda_n!1_{\lambda_1\mo_1\oplus\cdots\oplus\lambda_n\mo_n}$
for $\sum_{i=1}^n\lambda_i=m, \lambda_i\geq 0,i=1,...,n.$ Therefore
$$\bar{\varphi}:U(\Lambda)_{\mo_1\cdots \mo_n}^{(m)}/U(\Lambda)_{\mo_1\cdots
\mo_n}^{(m-1)}\rightarrow R(\Lambda)_{\mo_1\cdots
\mo_n}^{(m)}/R(\Lambda)_{\mo_1\cdots \mo_n}^{(m-1)}$$ is an
embedding for all $m\in\bbn.$ By a usual method of filtrated ring,
we see that
$$\varphi\mid_{U(\Lambda)_{\mo_1\cdots\mo_n}}:
U(\Lambda)_{\mo_1\cdots\mo_n}\rightarrow
R(\Lambda)_{\mo_1\cdots\mo_n}$$ is an embedding. It follows from
the arbitrary choice of $\mo_1,\cdots,\mo_n$ that
$\varphi:U(\Lambda)\rightarrow R(\Lambda)$ is an embedding.
\end{proof}

We define the subalgebra $R'(\Lambda)$ of $R(\Lambda)$ to be
generated by the elements
$\lambda_1!\cdots\lambda_n!1_{\lambda_1\mo_1\oplus\cdots\oplus\lambda_n\mo_n},$
where $\lambda_i\in\bbn$ for $i=1,\cdots,n$ and $\mo_i,
i=1,\cdots,n,$ are indecomposable constructible sets and disjoint to
each other.
\begin{Thm}\label{Rie} The canonical embedding $\varphi:U(\Lambda)\rightarrow
R(\Lambda)$ induces the isomorphism $\varphi:U(\Lambda)\rightarrow
R'(\Lambda).$ Therefore $U(\Lambda)\otimes_{\bbz}\bbq\cong
R(\Lambda)\otimes_{\bbz}\bbq.$
\end{Thm}

\begin{proof}
By Corollary \ref{coro3.7},
$$\varphi(1_{\mo_1}^{\ast\lambda_1}\ast\cdots\ast
1_{\mo_n}^{\ast\lambda_n})=1_{\mo_1}^{\lambda_1}\bullet\cdots\bullet
1_{\mo_n}^{\lambda_n}$$ $$=\lambda_1!\cdots\lambda_n!
1_{\lambda_1\mo_1\oplus\cdots\oplus\lambda_n\mo_n}+\sum_{\gamma(\mp_j)<m}c_{\mp_j}1_{\mp_j}$$
where $m=\sum_{i=1}^n\lambda_i, \ c_{\mp_j}\in\bbz$ and we may
assume that
$\mp_j=\lambda'_1\mo'_1\oplus\cdots\oplus\lambda'_t\mo'_t$,
$\sum_{i=1}^t\lambda'_i<m$ such that $\mo'_1,\cdots,\mo'_t$ are
indecomposable constructible sets and disjoint to each other. By
Lemma 2.2 in \cite{R} we know that $c_{\mp_j}$ can be divided by
$\lambda'_1!\cdots\lambda'_t!.$ By induction,
$c_{\mp_j}1_{\mp_j}\in\im(\varphi).$ Thus
$\lambda_1!\cdots\lambda_n!1_{\lambda_1\mo_1\oplus\cdots\oplus\lambda_n\mo_n}\in\im(\varphi).$
Obviously
\{$\lambda_1!\cdots\lambda_n!1_{\lambda_1\mo_1\oplus\cdots\oplus\lambda_n\mo_n}$\}
generate the whole $R(\Lambda)$ over $\bbq.$ So we have
$U(\Lambda)\otimes_{\bbz}\bbq\cong R(\Lambda)\otimes_{\bbz}\bbq.$
\end{proof}

\section{The Degenerated Form of Green's Formula}
\nd{\bf{5.1}} Let $f:Y\longrightarrow X$ be a homomorphism of
$\Lambda-$modules and $M=\im(f)$. For any $L\in \mod\Lambda$ with
$\udim L=\udim X+\udim Y-\udim M$, let $\mf(f;L)$ be the set of
all pairs $(c,d)$ where $c:L\longrightarrow X$ is an epimorphism,
$d:Y\longrightarrow L$ a monomorphism and $f=cd$, then
$\mathcal{F}(f;L)$ is a constructible subset of
$\hom_{\Lambda}(L,X)\times \hom_{\Lambda}(Y,L).$ Moreover we take
the $I$-graded space $V$ such that $\udim V=\udim X+\udim Y-\udim
M.$ A $\Lambda$-module structure on $V$ is $(V,\sigma)$ where
$\sigma:\Lambda\ra\ed_{\bbc}V$ is an algebra homomorphism. We
define $\mathcal{F}(f)=\bigcup_{\sigma:\Lambda\ra
{End_{\mathbb{C}}}V}\mathcal{F}(f;(V,\sigma))$, then
$\mathcal{F}(f)$ is a constructible subset of
$\hom_{\mathbb{C}}(V,X)\times \hom_{\mathbb{C}}(Y,V).$  Let
$$\gamma_1:0\longrightarrow \Ker f\longrightarrow Y\longrightarrow M\longrightarrow 0\in
 \ext_{\Lambda}^{1}(M,\Ker f)$$
$$\gamma_2:0\longrightarrow M\longrightarrow X\longrightarrow \Coker f\longrightarrow 0\in
 \ext_{\Lambda}^{1}(\Coker f,M).$$
Then naturally the equivalence class of
$\gamma=\gamma_2\gamma_1:0\longrightarrow \Ker f\longrightarrow
Y\longrightarrow X\longrightarrow \Coker f\longrightarrow 0$
belongs to $
 \ext_{\Lambda}^{2}(\Coker f,\Ker f).$
\begin{Lemma}\label{1} We have that
$\mathcal{F}(f)\not=\varnothing$ if and only if $\gamma=0.$
\end {Lemma}
\begin{proof}If $\gamma=0$, by the same discussion as in \cite{Ri3} we can
prove $\mathcal{F}(f)\not=\varnothing$. Conversely, if
$\mathcal{F}(f)\not=\varnothing$, then there exists $(L,c,d)$ such
that we have the following commutative diagram:
\[\begin{CD}
\gamma_0:~0 @>>> Y @>d>> L @>>> \Coker f @>>> 0\\
@. @VVV @VcVV @| @.\\
\gamma_2:~0 @>>> M @>>> X @>>> \Coker f @>>> 0
\end{CD}
\]
This means that for the following long exact sequence
$$\cdots\longrightarrow
\ext_{\Lambda}^{1}(\Coker f,\Ker f)\longrightarrow
\ext_{\Lambda}^{1}(\Coker f,Y)\stackrel{f}{\longrightarrow}$$
$$
\ext_{\Lambda}^{1}(\Coker f,M)\stackrel{g}{\longrightarrow}
\ext_{\Lambda}^{2}(\Coker f,\Ker f)$$  we have
$f(\gamma_0)=\gamma_2.$ Hence
$\gamma=\gamma_2\gamma_1=g(\gamma_2)=gf(\gamma_0)=0.$
\end{proof}

If $f$  factors through  $L$, we call $\gamma=0$ induced by $L.$

\bigskip
\nd{\bf{5.2 The description of $\ext_{\Lambda}^{1}(B,A)_{X}$}}. For
$A,B,X\in \mod \Lambda$, we define $\ext_{\Lambda}^{1}(B,A)_{X}$ to
be the subset of $\ext_{\Lambda}^{1}(B,A)$ with the middle term
isomorphic to $X.$ By the vector space structure of
$\ext_{\Lambda}^{1}(B,A)$, we can view $\ext_{\Lambda}^{1}(B,A)$ as
affine space $\mathbb{C}^{n}$ (
$n=\dim_{\mathbb{C}}\ext_{\Lambda}^{1}(B,A)$). We can prove
$\ext_{\Lambda}^{1}(B,A)_{X}$ is a constructible subset. Of course,
if $\udim X\not=\udim A+\udim B$, then
$\ext_{\Lambda}^{1}(B,A)_{X}=\varnothing$. Let $\alpha:\Lambda\ra
\ed_{\mathbb{C}}A$ and $\beta:\Lambda\ra \ed_{\mathbb{C}}B$ be the
algebra homomorphisms which describe the $\Lambda$-module structures
of $A$ and $B$, respectively. Set
$$Z(B,A)=\{\delta:\Lambda\ra
\hom_{\mathbb{C}}{(B,A)}:\delta(\lambda\mu)=\alpha(\lambda)\delta(\mu)+\delta(\lambda)\beta(\mu)\
\mbox{for}\
 \lambda,\mu\in \Lambda\}$$
$$T(B,A)=\{\delta:\Lambda\ra
\hom_{\mathbb{C}}{(B,A)}: \mbox{there exists}\ \eta\in
\hom_{\mathbb{C}}{(B,A)}\ \mbox{such that}$$$$\ \
\delta(\lambda)=\alpha(\lambda)\eta-\eta\beta(\lambda)\ \mbox{for
all} \ \lambda\in \Lambda\}$$ It is clear that both $Z(B,A)$ and
$T(B,A)$ are constructible subsets of \\
$\hom_{\mathbb{C}}(\Lambda,\hom_{\mathbb{C}}(B,A))$. There is a
natural morphism of varieties  (see \cite{R}) $\pi:Z(B,A)\ra
\ext_{\Lambda}^{1}(B,A).$ For $\delta\in Z(B,A)$, $\pi(\delta)$ is
the class of the following exact sequence in
$\ext_{\Lambda}^{1}(B,A)$:
$$0\longrightarrow A\stackrel{{1 \choose 0}}{\longrightarrow} (A\oplus B)_{\delta}\stackrel{(0\ \ 1)}
{\longrightarrow} B\longrightarrow0$$ where $(A\oplus B)_{\delta}$
is the $\mathbb{C}$-vector space $A\oplus B$ together with the
algebra homomorphism $\Lambda\ra \ed_{\mathbb{C}}(A\oplus B)$ which
sends $\lambda$ to \[
\left(\begin{array}{cc} \alpha(\lambda) & \delta(\lambda)\\
0 & \beta(\lambda) \end{array}\right)
\]
Any fibre of $\pi$ is homeomorphic to $T(B,A)$. Hence, if we set
$Z(B,A)_{X}=\{\delta\in Z(B,A)|(A\oplus B)_{\delta}\cong X\}$ and
$T(B,A)_{X}=T(B,A)\cap Z(B,A)_{X}$, then $\pi$ induces
$\pi_{X}:Z(B,A)_{X}\ra \ext_{\Lambda}^{1}(B,A)_{X}$ with the fibre
$T(B,A)_{X}$.
\begin{Lemma}\label{2}
If $X$ is not isomorphic to $A\oplus B,$ then
$\chi(\mathrm{Ext}_{\Lambda}^{1}(B,A)_{X})=0$.
\end{Lemma}
\begin{proof} Since $X\ncong A\oplus B$, $T(B,A)_{X}=\varnothing$ . This
shows $\chi(Z(B,A)_{X})=\chi(Ext_{\Lambda}^{1}(B,A)_{X})$. We
claim that $\chi(Z(B,A)_{X})=0$. If $\delta=0$, then $(A\oplus
B)_{\delta}\cong A\oplus B$ as $\Lambda-$module, so $0\notin
Z(B,A)_{X}$. There exists a free action of $\mathbb{C}^{*}$ on
$Z(B,A)_{X}$ by $t\cdot\delta=t\delta$ for $t\in \mathbb{C}^{*},$
whose orbit space is a constructible subset of projective space.
So
$\chi(Z(B,A)_{X})=\chi(\mathbb{C}^{*})\chi(Z(B,A)_{X}/\mathbb{C}^{*})=0$
by Proposition \ref{Euler}.
\end{proof}
\begin{Remark}
The above free action of $\mathbb{C}^{*}$ on $Z(B,A)_{X}$ induces
the following free action of $\mathbb{C}^{*}$ on
$\mathrm{Ext}^{1}_{\Lambda}(B,A)_X$(see 4.3 in \cite{R}): for
$t\in \mathbb{C}^{*}$ and an extension
$$
\varepsilon: 0\rightarrow A
\xrightarrow{a}Y\xrightarrow{b}B\rightarrow 0
$$
$$
t\cdot\varepsilon: 0\rightarrow A
\xrightarrow{t^{-1}a}Y\xrightarrow{b}B\rightarrow 0
$$
\end{Remark}

\bigskip
\nd{\bf{5.3}}  Let $X$ be a variety under the action of $G_{\ud}$
and $f$ be a $G_{\ud}$-invariant constructible function over $X.$
It is then natural to define $\int_{x\in X} f(x):=\sum_{m\in
\mathbb{C}}m\chi(f^{-1}(m)).$ Let $\overline{\mo}$ be the quotient
space of the constructible subset $\mo$ of $X$ under the action of
$G_{\ud}.$ Define
$$
\int_{x\in \overline{X}}f(x):=\sum_{m\in
\mathbb{C}}m\chi(\overline{f^{-1}(m)})
$$
where $\chi(\overline{f^{-1}(m)})$ is defined in Section 2.

Let $\mp$ be the set of isomorphism classes of $\Lambda$-module
and $V_{\alpha}$ be a representative in $\alpha$ for
$\alpha\in\mp.$ The dimension vector of $V_{\alpha}$ is denoted by
$\underline{\alpha}.$ Let $\mo$ be a $G_{\ud}$-invariant
constructible subset of $\bbe_{\ud}(Q,R)$, we set
$h^{\alpha\beta}_{\mo}=\chi(\ext_{\Lambda}^{1}(V_{\alpha},V_{\beta})_{\mo}),$
$g^{\lambda}_{\mo_{\beta}\mo_{\alpha}}=\chi(V(\mo_{\alpha},\mo_{\beta};V_{\lambda}))$
,
$h_{\lambda}^{\alpha\beta}=\chi(\ext_{\Lambda}^{1}(V_{\alpha},V_{\beta})_{V_{\lambda}})$
and
$g^{\lambda}_{\beta\alpha}=\chi(V(V_{\alpha},V_{\beta};V_{\lambda})).$

We fix $\mathcal{O}_{\alpha},\mathcal{O}_{\beta},\alpha^{\prime},$
$\beta^{\prime}$ where $\alpha^{\prime},\beta^{\prime}\in
\mathcal{P}$ and $\mathcal{O}_{\alpha},\mathcal{O}_{\beta}$ are
$G_{\underline{d}}$-invariant  constructible sets.
\begin{Lemma}\label{condition}
Given the following commutative diagram with exact rows and columns:
\[\begin{CD}
T @>u_T>> B @>q_T>> B/T\\
@Ve_1VV @. @Ve_2VV\\
B' @. @. A'\\
@Ve_3VV @. @Ve_4VV\\
S @>u_S>> A @>q_S>>A/S
\end{CD}
\]

\nd where $A^{\prime}\in \alpha^{\prime},B^{\prime}\in
\beta^{\prime}$, $A\in \mathcal{O}_{\alpha},B\in
\mathcal{O}_{\beta}$, $u_T, u_S$ are  canonical inclusions, and
 $q_T, q_S$ are  canonical projections.  We define the
pushout $$Y=B^{\prime}\sqcup_{T}B=B\oplus B'/\{e_1(t)\oplus
u_T(t)\mid t\in T\}$$ and the pullback
$$X=A^{\prime}\sqcap_{A'/S}A=\{(a\oplus a')\in A\oplus A'\mid e_4(a')=q_S(a)\}.$$ Then there exists the following exact
sequence:$$0\longrightarrow T\longrightarrow
Y\stackrel{f}{\longrightarrow} X\longrightarrow A/S\longrightarrow
0$$
\end{Lemma}
\begin{proof} See Section 6 in \cite{Ri3}.
\end{proof}
\bigskip

We denote by $\varepsilon(A/S,S,B/T,T)$ the equivalence class of
the above exact sequence in $\ext_{\Lambda}^{2}(A/S,T).$  By Lemma
5.1, if $\varepsilon(A/S,S,B/T,T)=0,$ then there exist $L$,
$d:Y\hookrightarrow L$ and $c:L\twoheadrightarrow X$ such that
$f=cd.$ We call $(S,T,e_1,e_2,e_3,e_4,c,d)$ is induced by $L.$

Let $\bbc^d$ and $\bbc^{d'}$ be two vector spaces of dimension $d$
and $d'$, respectively. Let $M_{d'\times d}$ be the set of all
matrices of size $d'\times d.$ Then $M_{d'\times d}=\hom(\bbc^d,
\bbc^{d'})$ and $M_{d'\times d}=\bigsqcup_{r}M_{d'\times d}(r)$
where $M_{d'\times d}(r)$ consists of all matrices of rank $r$.
For any $A=(a_{ij})\in M_{d'\times d}(r),$ let us denote the
$r\times r$ submatrix of $A$ formed by the rows $1\leq i_1<\cdots
<i_r\leq d'$ and the columns $1\leq j_1<\cdots <j_r\leq d$ by
$\bigtriangleup_{(i_1,\cdots, i_r; j_1,\cdots, j_r)}(A)$. For
every pair of multi-indices $I=\{i_1,\cdots, i_r\}\subseteq \{1,
\cdots, d'\}$ and $J=\{j_1,\cdots, j_r\}\subseteq \{1, \cdots,
d\},$ we define $M_{d'\times d}(r,I,J)$ to be the subset of
$M_{d'\times d}(r)$ consisting of the matrices $A$ which satisfy
$A\notin M_{d'\times d}(I', J')$ for any $I'<I$ or $I'=I, J'<J$
and $\mathrm{det}\ \bigtriangleup_{(i_1,\cdots, i_r; j_1,\cdots,
j_r)}(A)\neq 0$. Here $I'<I$ is the common lexicographic order. We
have a finite stratification of $M_{d'\times d}(r)$, i.e.,
$$
M_{d'\times d}(r)=\bigsqcup_{(I,J)}M_{d'\times d}(r,I,J).
$$
In particular, if $d<d'$ and $r=d$, this gives a finite
stratification of the Grassmannian $\mathrm{Gr}_{d}(\bbc^{d'})$
consisting of all $d$-dimensional subspaces of $\bbc^{d'}.$
Indeed, for any $I=\{i_1,\cdots, i_d\}\subset \{1, \cdots, d'\}$
and $J=\{1,\cdots, d\}$, let $M^g_{d'\times d}(I)$ be the subset
of $M_{d'\times d}(d, I,J)$ consisting of the matrices $A$
satisfying that $\bigtriangleup_{(I;J)}(A)$ are identity matrices.
Then there is a finite stratification
$$
\mathrm{Gr}_{d}(\bbc^{d'})=\bigsqcup_{I}M^g_{d'\times d}(I).
$$
For any $A\in M_{d'\times d}(d, I, J)$, we substitute the identity
matrix for the submatrix $\bigtriangleup_{(I;J)}(A)$ and then $A$
corresponds to a unqiue matrix $A'\in M^g_{d'\times d}(I)$.

For every pair of multi-indices $I=\{i_1,\cdots, i_r\}$ and
$J=\{j_1,\cdots, j_r\},$ we will define the following morphism of
varieties:
 $$\Upsilon_{(r, I, J)}^{1}: M_{d'\times d}(r, I,
J)\rightarrow M_{d\times (d-r)}(d-r),$$
$$\hspace{-0.8cm}\Upsilon_{(r, I,
J)}^{2}: M_{d'\times d}(r, I, J)\rightarrow
\mathrm{Gr}_{d-r}(\bbc^{d}),$$
$$\hspace{0.3cm}\Omega_{(r, I, J)}^{1}: M_{d'\times
d}(r, I, J)\rightarrow M_{(d'-r)\times d'}(d'-r), $$ and
$$\hspace{-1.2cm}\Omega_{(r, I, J)}^2: M_{d'\times
d}(r, I, J)\rightarrow \mathrm{Gr}_{r}(\bbc^{d'}).$$ Let
$P_{ij}(k)$ be the elementary matrix of size $k\times k$
transposing the $i$-th row and the $j$-th row. Set
$P_{I}(d')=P_{r,i_r}(d')\cdots P_{1, i_1}(d')$ and
$P_{J}(d)=P_{r,j_r}(d)\cdots P_{1, j_1}(d)$. Then we have
$$P_{I}(d')AP_J(d)\in M_{d'\times d}(r,(1,\cdots,r)(1,\cdots,r))$$ for
any matrix $A\in
M_{d'\times d}(r, I, J)$. The matrix $P_{I}(d')AP_J(d)$ has the form $\left(%
\begin{array}{cc}
  A_1 & A_2 \\
  A_3 & A_4 \\
\end{array}%
\right)$ with  an invertible $r\times r$ matrix $A_1$ and
$A_4=A_3A^{-1}_1A_2=A_2A^{-1}_1A_3.$ The matrix $P_J(d)\left(%
\begin{array}{c}
  -A_1^{-1}A_2 \\
  I_{d-r} \\
\end{array}%
\right)$ determines the solution space $\{x\in \bbc^{d}\mid
Ax=0\}$. The matrix $(-A_3A^{-1}_1, I_{d'-r})P_I(d')$ determines
the solution space $\{x\in \bbc^{d}\mid xA=0\}$.
We define $$\Upsilon^1_{(r, I, J)}(A)=P_J(d)\left(%
\begin{array}{c}
  -A_1^{-1}A_2 \\
  I_{d-r} \\
\end{array}%
\right), $$ $$\Omega^1_{(r, I, J)}(A)=(-A_3A^{-1}_1,
I_{d'-r})P_I(d').$$

Assume that $P_J(d)\left(%
\begin{array}{c}
  -A_1^{-1}A_2 \\
  I_{d-r} \\
\end{array}%
\right)\in M_{d\times (d-r)}(d-r, I', J')$ for some $I'\subset
\{1,\cdots, d\}$ and $J'=(1,\cdots, d-r)$.  Then we define
$\Upsilon^2_{(r, I, J)}(A)$ to be the
unique matrix in $M^g_{d\times (d-r)}(I')$ which the matrix $P_J(d)\left(%
\begin{array}{c}
  -A_1^{-1}A_2 \\
  I_{d-r} \\
\end{array}%
\right)$ corresponds to. Similarly, we define $\Omega^2_{(r, I,
J)}(A)$ to be the unique matrix in $M^g_{d'\times r}(I)$ which the
submatrix $\bigtriangleup_{(1,\cdots, d'; j_1,\cdots, j_r})(A)$ of
$A$ corresponds to. Hence, for any $A\in M_{d'\times d}(r, I, J)$,
we have a long exact sequence of $\bbc$-spaces
$$
\xymatrix{0\ar[r]&\bbc^{d-r}\ar[rr]^{\Upsilon^1_{(r, I,
J)}(A)}&&\bbc^{d}\ar[r]^{A}&\bbc^{d'}\ar[rr]^{\Omega^1_{(r, I,
J)}(A)}&&\bbc^{d'-r}\ar[r]&0}.
$$

We can replace the vector spaces by $\Lambda$-modules. The
discussion is similar.  It gives a intrinsic description of kernel
and cokernel. We have the following lemma.
\begin{Lemma}\label{kernel}
Let $X\in \bbe_{d_1}(\Lambda)$ and $Y\in \bbe_{\ud_2}(\Lambda).$
Then there are the following maps whose restrictions to the strata
are morphisms of varieties
$$
\mathcal{K}:\mathrm{Hom}_{\Lambda}(X,Y)\rightarrow \bigcup_{\ue}
Gr_{\ue}(X) \quad \mbox{ and }\quad G: \bigcup_{\ue}
Gr_{\ue}(X)\rightarrow \bigcup_{\ue<\ud_1}\bbe_{\ue}(\Lambda)
$$
where the first morphism is taking the kernel and the second
morphism is the substitution of underlying space. Dually, there
exists a map whose restriction to the strata are  morphisms of
varieties
$$
\mathcal{C}:\mathrm{Hom}_{\Lambda}(X,Y)\rightarrow
\bigcup_{\ue'<\ud_2} \bbe_{\ue'}(\Lambda).
$$
defined by taking cokernel.
\end{Lemma}

Given  $\Lambda$-modules $A,B, A', B'$  and $L,$ we define
$\mathcal{Q}(L)$ to be the set of all quadruples $(a,b,a',b')$
which satisfy the following diagrams
\[\begin{CD}
@. B @.\\
@. @VVbV @.\\
B' @>b'>> L @>a'>> A'\\
@. @VVaV @.\\
@. A @.
\end{CD}
\] where both rows and columns are short exact sequences.

Define $\mathcal{O}(L)$ to be the set of all octets
$(S,T,e_1,e_2,e_3,e_4, c,d)$ which satisfy the following diagram
$$
\xymatrix{ & 0 \ar[d] &&&&  0\ar[d] && \\
0 \ar[r] & T \ar[rr]^-{u_T} \ar[dd]^-{e_1} && B\ar[rr]^-{q_T} \ar@{.>}[dl]\ar@{.>}[dd]&& B/T \ar[r] \ar[dd]^-{e_3}& 0\\
&&Y\ar@{.>}[dr]^-{d}&&&&\\
& B' \ar@{.>}[ur]\ar[dd]^-{e_2} \ar@{.>}[rr]&& L\ar@{.>}[rr]\ar@{.>}[dd]\ar@{.>}[dr]^-{c}&& A' \ar[dd]^-{e_4} &\\
&&&&X\ar@{.>}[ur]\ar@{.>}[dl]&&\\
0\ar[r] &  S \ar[d] \ar[rr]^-{u_S} && A \ar[rr]^-{q_S} && A/S \ar[r] \ar[d] & 0\\
& 0 && && 0 &}
$$ Here both rows and columns are short exact sequences, $S$ and $T$ are submodules of $A$ and $B$, respectively and
$u_T, u_S$ are the canonical inclusions, $q_T, q_S$ are the
canonical projections and $X,Y$ are as in Lemma \ref{condition}.
Moreover, $\varepsilon(A/S,S,B/T,T)=0,$ and any octet in
$\mathcal{O}(L)$ is induced by $L.$

We define a map $\eta:\mathcal{Q}(L)\longrightarrow
\mathcal{O}(L)$ as follows. For any $(a,b,a',b')\in
\mathcal{Q}(L),$ we set
$$S=ab'(B'),T=b^{-1}b'(B),
e_1=b'^{-1}b, e_2=a'b(q_T)^{-1}, e_3=ab', e_4=q_Saa'^{-1}$$ and
$d=(b,b')\mid_{(B\oplus B')/T}, c=(a,-a')$. Then, there exists an
exact sequence:
$$
0\longrightarrow T\longrightarrow (B\oplus
B')/T\stackrel{f}{\longrightarrow}\Ker(q_S,e_4)\longrightarrow
A/S\longrightarrow 0
$$ and $f=cd.$ This shows $\mathcal{F}(f)\not=\varnothing.$ By Lemma 5.1,
$\varepsilon(A/S,S,B/T,T)=0$. Hence, we define
$$
\eta((a,b,a',b'))=(S,T,e_1,e_2,e_3,e_4,c,d).
$$
By Lemma \ref{kernel}, $\eta$ is a morphism of varieties.

Conversely, let $i_B: B\rightarrow B\oplus B'$ and $i_{B'}:
B'\rightarrow B\oplus B'$ be the canonical embeddings and let
$p_B: B\oplus B'\rightarrow B$ and $p_{B'}: B\oplus B'\rightarrow
B'$ be the canonical projections,  we have a map
$\zeta:\mathcal{O}(L)\longrightarrow \mathcal{Q}(L)$ as follows.
Consider the morphisms of $\Lambda$-modules:
$$
\xymatrix{B\ar[r]^-{i_B}& B\oplus B'\ar[r]^-{\pi}& (B\oplus
B')/T\ar[r]^-{d}& L}
$$
and
$$
\xymatrix{B'\ar[r]^-{i_{B'}}& B\oplus B'\ar[r]^-{\pi}& (B\oplus
B')/T\ar[r]^-{d}& L}.
$$
Define
$$
b=d\circ\pi\circ i_B, \quad b'=d\circ\pi\circ i_{B'}
$$

Also, we consider the morphisms
$$
\xymatrix{L\ar[r]^-{c}& \mathrm{ker}(q_S, e_4)\ar[r]^-{i}& A\oplus
A'\ar[r]^-{p_A}& A}
$$
and
$$
\xymatrix{L\ar[r]^-{c}& \mathrm{ker}(q_S, e_4)\ar[r]^-{i}& A\oplus
A'\ar[r]^-{p_{A'}}& A'}.
$$
Define
$$
a=p_A\circ i\circ c, \quad a'=p_{A'}\circ i\circ c.
$$
Hence, we define
$$
\xi((S,T,e_1,e_2,e_3,e_4,c,d))=(a,b,a',b').
$$
By Lemma \ref{kernel}, $\xi$ is a morphism of varieties. Now we
have the following proposition which can be viewed as the
geometrization of the bijections in \cite{G} and \cite{Ri3}.
\begin{Prop}
There exists a  map from $Q(L)$ to $\mo(L)$ whose restriction to
the strata are homeomorphisms between varieties.
\end{Prop}

Our main theorem is the following formula which can be viewed as the
degenerated form of Green's formula.

\begin{Thm}\label{greentheorem}
Let
$\mathcal{O}_{\alpha},\mathcal{O}_{\beta},\alpha^{\prime},\beta^{\prime}$
as above, then
$$g_{\mathcal{O}_{\alpha}\mathcal{O}_{\beta}}^{\alpha'\oplus \beta'}
=\int_{(\rho,\sigma,\sigma^{\prime},\tau)\in \overline{\bbe},
\rho\oplus \sigma\in \overline{\mo}_{\alpha}, \sigma'\oplus
\tau\in \overline{\mo}_{\beta}}
g_{\sigma\tau}^{\beta^{\prime}}g_{\rho\sigma^{\prime}}^{\alpha^{\prime}}$$
where $\bbe=\bbe_{\underline{\rho}}(\Lambda)\times
\bbe_{\underline{\sigma}}(\Lambda)\times
\bbe_{\underline{\sigma'}}(\Lambda)\times
\bbe_{\underline{\tau}}(\Lambda)$.
\end{Thm}
\begin{proof}
It is enough to prove that for fixed $\Lambda$-modules
$V_{\alpha},V_{\beta},V_{\alpha'},V_{\beta'},$ the following
identity holds
$$\int_{\alpha\in \mo_{\alpha}, \beta\in\mo_{\beta}}g_{V_{\alpha}V_{\beta}}^{V_{\alpha'\oplus\beta'}}
=\int_{(\rho,\sigma,\sigma^{\prime},\tau)\in \overline{\bbe},
\rho\oplus \sigma\in \overline{\mo}_{\alpha}, \sigma'\oplus
\tau\in \overline{\mo}_{\beta}}
g_{V_{\sigma}V_{\tau}}^{V_{\beta^{\prime}}}g_{V_{\rho}V_{\sigma^{\prime}}}^{V_{\alpha^{\prime}}}.$$
There is a natural embedding morphism:
$$
\bigcup_{\rho,\sigma,\sigma',\tau;\rho\oplus
\sigma=\alpha,\sigma'\oplus
\tau=\beta}V(V_{\sigma'},V_{\rho};V_{\alpha'})\times
V(V_{\tau},V_{\sigma};V_{\beta'})\xrightarrow{i}
V(V_{\beta},V_{\alpha};V_{\alpha'}\oplus V_{\beta'}).
$$
Consider the action of $\bbc^{*}$ on
$V(V_{\beta},V_{\alpha};V_{\alpha'}\oplus V_{\beta'})$ as in the
proof of Proposition \ref{action}. Its fixed point set is just
$\mbox{Im}i.$ The proof of the above identity follows from the
fact that $V(V_{\beta},V_{\alpha};V_{\alpha'}\oplus V_{\beta'})$
has the same Euler characteristic as its fixed point set under the
action of $\bbc^{*}.$
\end{proof}

\begin{Remark}\label{extended}
\nd (1) One can prove the following version of Green's formula in
case the base field $k$ is a finite field with the cardinality
$q$:
$$a_{\alpha}a_{\beta}a_{\alpha^{'}}a_{\beta^{'}}
\sum_{\lambda}g_{\alpha\beta}^{\lambda}g_{\alpha^{'}\beta^{'}}^{\lambda}a_{\lambda}^{-1}$$
$$=
\sum_{\rho,\sigma,\sigma^{'},\tau,\varepsilon(\rho,\sigma,\sigma^{'},\tau)=0}
\frac{|\mathrm{Ext}^{1}(V_{\rho},V_{\tau})|}{|\mathrm{Hom}(V_{\rho},V_{\tau})|}
g_{\rho\sigma}^{\alpha}g_{\rho\sigma^{'}}^{\alpha^{'}}g_{\sigma^{'}\tau}^{\beta}g_{\sigma\tau}^{\beta^{'}}
a_{\rho}a_{\sigma}a_{\sigma^{'}}a_{\tau}$$ which extends Green's
formula (see \cite{G} and \cite{Ri3}).  In particular, when the
category of finite $\Lambda$-modules is hereditary,
$\varepsilon(\rho,\sigma,\sigma^{'},\tau)$ always vanishes so that
the above form is reduced to the original Green's formula. Note
that the original Green's formula holds for all finite
$\Lambda$-modules if and only if the category of finite
$\Lambda$-modules is hereditary (see \cite{Ri3} and \cite{ZW}).
\\
\nd (2) Let
$h_{\lambda}^{\alpha'\beta'}=|\mathrm{Ext}^1(V_{\alpha'},
V_{\beta'})_{\lambda}|.$ Then the above extended form of Green's
formula can be rewritten as
$$
\sum_{\lambda}g_{\alpha\beta}^{\lambda}h_{\lambda}^{\alpha'\beta'}$$
$$=
\sum_{\rho,\sigma,\sigma^{'},\tau,\varepsilon(\rho,\sigma,\sigma^{'},\tau)=0}
\frac{|\mathrm{Ext}^{1}(V_{\rho},V_{\tau})|\cdot|\mathrm{Hom}(V_{\alpha'},V_{\beta'})|}{|\mathrm{Hom}(V_{\rho},V_{\tau})|\cdot|\mathrm{Hom}(V_{\rho},V_{\sigma})|\cdot|\mathrm{Hom}(V_{\sigma'},V_{\tau})|}
g_{\rho\sigma^{'}}^{\alpha^{'}}g_{\sigma\tau}^{\beta^{'}}h_{\alpha}^{\rho\sigma}h_{\beta}^{\sigma'\tau}.$$
As Lemma \ref{2}, we know $h_{\lambda}^{\alpha'\beta'}\equiv
0\quad
 \mathrm{mod}(q-1)$ unless $\lambda=\alpha'\oplus\beta'.$ Hence,
 the above reformulation implies
 $$
g_{\alpha\beta}^{\alpha'\oplus\beta'}\equiv
\sum_{\rho,\sigma,\sigma^{\prime},\tau; \rho\oplus \sigma= \alpha,
\sigma'\oplus \tau=\beta}
g_{\sigma\tau}^{\beta^{\prime}}g_{\rho\sigma^{\prime}}^{\alpha^{\prime}}\quad
 \mathrm{mod}(q-1).
 $$
which inspires Theorem \ref{greentheorem}.
\\
\nd (3) If $\Lambda=\bbc Q$ for some acyclic quiver $Q$, let $X_M$
be the image of the Caldero-Chapton map (\cite{CC}) for any
$\Lambda$-module $M$. Then Theorem \ref{greentheorem} is
equivalent to the fact $X_MX_N=X_{M\oplus N}$ for any
$\Lambda$-modules $M$ and $N$ (\cite{CK}).
\end{Remark}
\section{Comultiplications over Universal Enveloping Algebras}

\nd As we know, $R(\Lambda)\otimes_{\bbz}\mathbb{Q}$ can be viewed
as the universal enveloping algebra of $L(\Lambda)\otimes_{\bbz}
\mathbb{Q}$. On the other hand, there is a unique comultiplication
over the universal enveloping algebra. We can realize this
comultiplication in our geometric frame of $R'(\Lambda)\ (\cong
U(\Lambda))$. Let $\delta:R'(\Lambda)\longrightarrow
R'(\Lambda)\otimes_{\mathbb{Z}}R'(\Lambda)$ be given by
$$\delta(\lambda_1!\cdots\lambda_n!1_{\mathcal{O}_{\lambda}})(A,B)=\chi(\ext_{\Lambda}^{1}(A,B)_{\mathcal{O}_{\lambda}})$$
for any  Krull-Schmidt constructible subset
$\mo_{\lambda}=\lambda_1\mo_1\oplus\cdots\oplus \lambda_n\mo_n.$
\begin{Prop} The operator
$\delta$ is well-defined i.e $\delta(1_{\mathcal{O}})$ is a
constructible function for any constructible set of stratified
Krull-Schmidt $\mathcal{O}$.
\end{Prop}
\begin{proof} We only need to consider the constructible set of
Krull-Schmidt. Let $\mo_{\lambda}=\lambda_1\mo_1\oplus\cdots\oplus
\lambda_n\mo_n$ be the Krull-Schmidt decomposition for the
constructible set $\mo_{\lambda}\subseteq \bbe_{\ud}(Q,R)$ and we
assume more that $\mo_i, 1\leq i\leq n,$ are disjoint to each
other. Then
$\lambda_1!\cdots\lambda_n!1_{\lambda_1\mo_1\oplus\cdots\oplus
\lambda_n\mo_n}\in R'(\Lambda).$

For a $G_{\underline{\lambda}}$-invariant constructible subset
$\mo_{\lambda}\subseteq \bbe_{\underline{\lambda}}(\Lambda),$ we
consider the following set
$$\ext_{\Lambda}^{1}(\mathbb{E}_{\underline{d}_2},
\mathbb{E}_{\underline{d}_1})_{\mathcal{O}_{\lambda}}=\{(A,B,\varepsilon)\mid
B\in \mathbb{E}_{\underline{d_1}},A\in
\mathbb{E}_{\underline{d_2}},\varepsilon\in
\ext^1_{\Lambda}(A,B)_{L}, L\in \mathcal{O}_{\lambda}\}$$ where
$\underline{\lambda}=\ud_1+\ud_2.$ It is a constructible subset of
$\ext^1_{\Lambda}(A,B)\times \bbe_{\ud_1}(\Lambda)\times
\bbe_{\ud_2}(\Lambda)$. We have the canonical projection
$$h: \ext_{\Lambda}^{1}(\mathbb{E}_{\underline{d}_2},
\mathbb{E}_{\underline{d}_1})_{\mathcal{O}_{\lambda}}{\longrightarrow}\mathbb{E}_{\underline{d}_1}\times
\mathbb{E}_{\underline{d}_2}.$$ For
$A\in\mathbb{E}_{\underline{d}_2},$
$B\in\mathbb{E}_{\underline{d}_1} ,\
h^{-1}(A,B)=\ext_{\Lambda}^{1}(A,B)_{\mathcal{O}_{\lambda}}$. By
Theorem \ref{Joyce}, we have
$h_{*}(1_{\ext_{\Lambda}^{1}(\mathbb{E}_{\underline{d}_2},
\mathbb{E}_{\underline{d}_1})_{\mathcal{O}_{\lambda}}})=\delta(1_{\mathcal{O}_{\lambda}})|_{\mathbb{E}_{\underline{d}_1}\times
\mathbb{E}_{\underline{d}_2}}$ is a constructible function. On the
other hand, there are finitely many
$\underline{d}_1,\underline{d}_2$ such that
$\underline{d}_1+\underline{d}_2=\underline{\lambda}$. Therefore,
$\delta(1_{\mathcal{O}_{\lambda}})$ is a constructible function.
Now we can write down:
$$
\delta(1_{\mathcal{O}_{\lambda}})=\sum_{i=1}^t
h_{\mo_{\lambda}}^{\alpha\beta}1_{\mo_i}
$$
where $\mo_i\in\bbe_{\ud_1}\times\bbe_{\ud_2}$ are constructible,
$\ud_1+\ud_2=\underline{\lambda}$, $(V_{\alpha}, V_{\beta})\in
\mo_i$ and
$h_{\mathcal{O}_{\lambda}}^{\alpha\beta}=\chi(\ext_{\Lambda}^{1}(V_{\alpha},V_{\beta})_{\mathcal{O}_{\lambda}})$.
By Lemma \ref{2}, we know that
$h_{\mathcal{O}_{\lambda}}^{\alpha\beta}=1$ if $V_{\alpha}\oplus
V_{\beta} \in \mo_{\lambda}$, otherwise $0$. We can write down the
formula explicitly:
$$
\delta(\lambda_1!\cdots\lambda_n!1_{\mo_{\lambda}})=\sum_{k_1,\cdots,k_n;
k_i\leq
\lambda_i}\lambda_1!\cdots\lambda_n!(1_{k_1\mo_1\oplus\cdots\oplus
k_n\mo_n}\otimes 1_{(\lambda_1-k_1)\mo_1\oplus\cdots\oplus
(\lambda_n-k_n)\mo_n}).
$$
The right side of the identity can be written as $$
\sum_{k_1,\cdots,k_n; k_i\leq \lambda_i}C^{k_1}_{\lambda_1}\cdots
C^{k_n}_{\lambda_n}\cdot$$ $$ k_1!\cdots
k_n!1_{k_1\mo_1\oplus\cdots\oplus k_n\mo_n}\otimes
(\lambda_1-k_1)!\cdots
(\lambda_n-k_n)!1_{(\lambda_1-k_1)\mo_1\oplus\cdots\oplus
(\lambda_n-k_n)\mo_n}. $$If $\mo_{\lambda}$ is indecomposable then
$\delta(1_{\mathcal{O}})=1_{\mathcal{O}}\otimes 1+1\otimes
1_{\mathcal{O}}$ by Lemma \ref{2}. The proof is finished.
\end{proof}
\begin{Thm}
The map $\delta:R(\Lambda)\longrightarrow R(\Lambda)\otimes
R(\Lambda)$ (also $\delta:R'(\Lambda)\longrightarrow
R'(\Lambda)\otimes R'(\Lambda)$) is an algebra homomorphism.
\end{Thm}
\begin{proof}Consider the following diagram:

\[
\begin{CD}
R(\Lambda)\otimes R(\Lambda) @>\mu>> R(\Lambda)\\
@V\delta\otimes \delta VV @VV\delta V\\
R(\Lambda)\otimes R(\Lambda)\otimes R(\Lambda)\otimes R(\Lambda)
@>\eta>> R(\Lambda)\otimes R(\Lambda)
\end{CD}
\]
where $\eta=(\mu\otimes\mu)\xi$:$$R(\Lambda)\otimes
R(\Lambda)\otimes R(\Lambda)\otimes
R(\Lambda)\stackrel{\xi}{\longrightarrow}R(\Lambda)\otimes
R(\Lambda)\otimes R(\Lambda)\otimes
R(\Lambda)\stackrel{\mu\otimes\mu}{\longrightarrow}R(\Lambda)\otimes
R(\Lambda)$$ given by $\xi(x\otimes y \otimes z\otimes w)=x\otimes
z \otimes y\otimes w.$ We need to show this diagram commutes. For
$1_{\mathcal{O}_{1}} \otimes 1_{\mathcal{O}_{2}}\in
R(\Lambda)\otimes R(\Lambda)$, using Lemma \ref{2}, we have
\begin{eqnarray}
  \delta\mu(1_{\mathcal{O}_{\alpha}} \otimes1_{\mathcal{O}_{\beta}}) &=& \delta(\sum_{\lambda}g_{\mathcal{O}_{\alpha}\mathcal{O}_{\beta}}^{\lambda}1_{\mathcal{O}_{\lambda})}=\sum_{\lambda}g_{\mathcal{O}_{\alpha}\mathcal{O}_{\beta}}^{\lambda}\sum_{\alpha^{'},\beta^{'}}
h_{\mathcal{O}_{\lambda}}^{\alpha^{'}\beta^{'}}1_{\mathcal{O}_{\alpha^{'}}}
\otimes 1_{\mathcal{O}_{\beta^{'}}} \nonumber\\
   &=& \sum_{\alpha^{'},\beta^{'}}(\sum_{\lambda}g_{\mathcal{O}_{\alpha}\mathcal{O}_{\beta}}^{\lambda})
h_{\mathcal{O}_{\lambda}}^{\alpha^{'}\beta^{'}}1_{\mathcal{O}_{\alpha^{'}}}
\otimes
1_{\mathcal{O}_{\beta^{'}}}=\sum_{\alpha^{'},\beta^{'}}g_{\mathcal{O}_{\alpha}\mathcal{O}_{\beta}}^{\alpha'\oplus\beta'}
1_{\mathcal{O}_{\alpha^{'}}} \otimes
1_{\mathcal{O}_{\beta^{'}}}\nonumber
\end{eqnarray}
and
\begin{eqnarray}
   && (\mu\otimes\mu)\xi(\delta\otimes\delta)(1_{\mathcal{O}_{\alpha}} \otimes1_{\mathcal{O}_{\beta}})=
(\mu\otimes\mu)\xi(\delta(1_{\mathcal{O}_{\alpha}})\otimes\delta(1_{\mathcal{O}_{\beta}}))\nonumber \\
  &=& (\mu\otimes\mu)\xi(\sum_{\rho,\sigma,\sigma^{'},\tau}h_{\mathcal{O}_{\alpha}}^{\rho\sigma}
h_{\mathcal{O}_{\beta}}^{\sigma^{'}\tau}1_{\mathcal{O}_{\rho}}
\otimes1_{\mathcal{O}_{\sigma}}\otimes1_{\mathcal{O}_{\sigma^{'}}}
\otimes1_{\mathcal{O}_{\tau}}) \nonumber\\
   &=& (\mu\otimes\mu)(\sum_{\rho,\sigma,\sigma^{'},\tau}h_{\mathcal{O}_{\alpha}}^{\rho\sigma}
h_{\mathcal{O}_{\beta}}^{\sigma^{'}\tau}1_{\mathcal{O}_{\rho}}
\otimes1_{\mathcal{O}_{\sigma^{'}}}\otimes1_{\mathcal{O}_{\sigma}}
\otimes1_{\mathcal{O}_{\tau}})\nonumber \\
   &=& \sum_{\rho,\sigma,\sigma^{'},\tau}h_{\mathcal{O}_{\alpha}}^{\rho\sigma}
h_{\mathcal{O}_{\beta}}^{\sigma^{'}\tau}(\sum_{\alpha^{'}}g_{\mathcal{O}_{\rho}\mathcal{O}_{\sigma^{'}}}^{\alpha^{'}}
1_{\mathcal{O}_{\alpha^{'}}})\otimes(\sum_{\beta^{'}}g_{\mathcal{O}_{\sigma}\mathcal{O}_{\tau}}^{\beta^{'}}
1_{\mathcal{O}_{\beta^{'}}}) \nonumber\\
   &=&\sum_{\alpha^{'},\beta^{'}}(\sum_{\rho,\sigma,\sigma^{'},\tau}h_{\mathcal{O}_{\alpha}}^{\rho\sigma}
h_{\mathcal{O}_{\beta}}^{\sigma^{'}\tau}g_{\mathcal{O}_{\rho}\mathcal{O}_{\sigma^{'}}}^{\alpha^{'}}
g_{\mathcal{O}_{\sigma}\mathcal{O}_{\tau}}^{\beta^{'}})1_{\mathcal{O}_{\alpha^{'}}}
\otimes 1_{\mathcal{O}_{\beta^{'}}}\nonumber\\
   &=&\sum_{\alpha^{'},\beta^{'}}(\int_{\rho,\sigma,\sigma^{'},\tau; \rho\oplus
\sigma\in \overline{\mo}_{\alpha}, \sigma'\oplus \tau\in
\overline{\mo}_{\beta}}g_{\rho\sigma^{'}}^{\alpha^{'}}
g_{\sigma\tau}^{\beta^{'}})1_{\mathcal{O}_{\alpha^{'}}} \otimes
1_{\mathcal{O}_{\beta^{'}}}\nonumber
\end{eqnarray}

Therefore, by Theorem \ref{greentheorem}, we have
$$
\delta\mu(1_{\mathcal{O}_{\alpha}}
\otimes1_{\mathcal{O}_{\beta}})=
(\mu\otimes\mu)\xi(\delta\otimes\delta)(1_{\mathcal{O}_{\alpha}}
\otimes1_{\mathcal{O}_{\beta}}).$$ We complete the proof of this
theorem.
\end{proof}

%\begin{titlepage}
%\def\refname{\hfil REFERENCES}

%\end{titlepage}

\end{document}